\documentclass{amsart}

\usepackage{xcolor,latexsym,amsfonts,amssymb,bbm,
amsmath,cite, amsthm, bbold, float}
\usepackage{graphicx}
\usepackage{hyperref}
\usepackage{url}
\usepackage{soul}
\usepackage{comment}
\usepackage{dsfont}
\usepackage[T1]{fontenc}
\usepackage{bbding}
\usepackage{pifont}
\usepackage{wasysym}

\title{Structural results for the Tree Builder Random Walk}
\author{J. Engl\"ander}
\address{Department of Mathematics, University of Colorado,
 Boulder, CO-80309, USA}
\email{janos.englander@colorado.edu}

\author{G. Iacobelli}
\address{ Mathematical Institute,  Federal University of Rio de Janeiro (UFRJ) -  Brazil} 
\email{giulio@im.ufrj.br} 

\author{G. Pete}
\address{HUN-REN Alfr\'ed R\'enyi Institute of Mathematics, Hungary, and Department of Stochastics, Institute of Mathematics, Budapest University of Technology and Economics, M\H{u}egyetem rkp.~3., Budapest, 1111, Hungary}
\email{gabor.pete@renyi.hu}

\author{R. Ribeiro}
\address{Department of Mathematics, University of Denver,
 CO-80210, USA.}
\email{Rodrigo.Ribeiro@du.edu} 

\newtheorem{theorem}{Theorem}[section]
\newtheorem{lemma}{Lemma}[section]
\newtheorem{corollary}{Corollary}[section]
\newtheorem{proposition}{Proposition}[section]
\newtheorem{definition}{Definition}[section]
\newtheorem{remark}{Remark}[section]
\newtheorem{example}{Example}[section]
\newtheorem{problem}{Problem}[section]

\newcommand{\mixc}{(M)}
\newcommand{\laws}{\mathcal{L}}

\newcommand{\Pd}{\mathbb{P}}
\newcommand{\Ed}{\mathbb{E}}

\newcommand{\name}{TBRW}

\newcommand{\degree}[2]{\mathrm{deg}_{#1}(#2)}

\newcommand{\com}[1]{\textcolor{red}{\texttt{giulio:} #1}}

\newcommand{\veps}{\varepsilon} 
\newcommand{\eps}{\varepsilon}

\newcommand{\myexp}{2(1-\gamma)+\delta}
\newcommand{\cL}{\mathcal{L}}

\def\R{\mathbb{R}}
\def\g{g}
\def\G{G}

\begin{document}

\begin{abstract}We study the Tree Builder Random Walk: a randomly growing tree, built by a walker as she is walking around the tree. Namely, at each time $n$, she adds a leaf to her current vertex with probability $p_n \asymp n^{-\gamma}$, $\gamma\in (2/3,1]$, then moves to a uniform random neighbor on the possibly modified tree. We show that the tree process at its growth times, after a random finite number of steps, can be coupled to be identical to the Barab\'asi-Albert preferential attachment tree model. 

Thus, our TBRW-model is a local dynamics giving rise to the BA-model. The coupling also implies that many properties known for the BA-model, such as diameter and degree distribution, can be directly transferred to our TBRW-model, extending previous results.
\end{abstract}

\maketitle

\section{Introduction}

In this work, we explore a class of stochastic processes that sits at the intersection of growing random graph models and random walks in dynamic environments. From one perspective, this model can be interpreted as a random graph generated by a random walk. From another, it can be seen as a random walk in a dynamic environment, with the feature that the environment's dynamics is heavily influenced by the walker's trajectory. Namely, we investigate an instance of the so-called Tree Builder Random Walks (\name) \cite{IRVZ22}. In this model, a random number of vertices are added to the walker's position at each step, resulting in a tree whose structure depends on the walker's trajectory.

An intriguing aspect of this dual perspective, characteristic to \name{}, is its automatic inheritance of natural questions associated with both random graphs and random walks. Therefore, questions related to random graphs such as diameter size, degree distribution, and clustering make sense in the context of the \name. Likewise, questions related to random walks like transience, recurrence, and ballisticity are also pertinent to the \name{}.
In this work, we examine the \name{} from the perspective of a random graph. We will focus on investigating the  properties of the graphs generated by its mechanisms.

 The link between random walks and (random) graphs is not a new phenomenon. Random walks have been employed as tools to estimate network properties; for instance, in \cite{BOY2018} random walks are used to estimate the number of vertices and edges of networks, or the PageRank algorithm uses stationary distributions to define the importance of vertices \cite{gleich2015pagerank}. A classical example of random trees generated by random walks is the Uniform Spanning Tree of any finite graph generated by Wilson's or the Aldous-Broder algorithm \cite[Chapter 4]{MR3616205}. Self-repelling or reinforced random walks, where the transition probabilities on a given graph are evolving based on the history of the random walk also have a rich literature; two beautiful examples are \cite{toth1998true} and \cite{angel2014localization}. Moreover, it has been observed that real-world networks, from biology through game theoretical economics to the internet, are often grown by agents moving and acting locally in the existing network; see \cite{gross2008adaptive} for an interdisciplinary literature review. In particular, the idea that scale-free (power-law degree distribution) graphs should be built locally by random walks, instead of the global Barab\'asi-Albert preferential attachment rule \cite{barabasi1999emergence} and its generalizations, has been proposed several times \cite{vazquez2003growing,saramaki2004scale,evans2005scale,amorim2016growing, ross2019random}, and a solution to this general problem is one of the main contributions of our present work. A concrete recent application of the \name{} is that random walks that add new vertices can be useful in chemistry, especially for describing ``chain walking'' in catalysts during polymer creation \cite{DS22}.

In Subsections~\ref{ss:model} and~\ref{ss:BA}, we formally introduce our model and the Barab\'asi-Albert random tree, then in Subsection~\ref{ss:app} we discuss their connections and further applications of the \name{}. In Subsection~\ref{ss:res} we describe our results and some of the many open problems that remain. We close this introductory section in Subsection~\ref{ss:ideas} with a sketch of the proof ideas.

\subsection{The \name{} model}\label{ss:model}

The \name{} is a stochastic process $\{(T_n, X_n)\}_{n\geq 0}$ where $T_n$ is a tree and $X_n$, a vertex of the tree, is the current position of a random walker. The \name{} model depends  on a sequence of probability laws~$\mathcal{L}:=\mathcal{L}_1,\mathcal{L}_2,...$, with each $\cL_n$ supported on nonnegative integers,  and an initial condition $(T_0, x_0)$, where 
\begin{itemize}
\item[(a)] $T_0$ is a finite, rooted tree with a self-loop attached at the root (to avoid periodicity); 
\item[(b)] $x_0$ is a vertex of $T_0$. 
\end{itemize}
 At every time $n$, a random $Z_n$ number of leaves, distributed according to  $\mathcal{L}_n$, independent of everything else, are added to $X_{n-1}$, the current position of the walker. Whenever $Z_n>0$ happens, the tree $T_{n-1}$ is modified immediately. After that, the random walk takes a step on the possibly modified tree to a uniformly chosen neighbor, before the next new leaves are drawn. 

The \name{}
has been recently introduced in \cite{IRVZ22} under the assumption of uniform ellipticity on the sequence of laws $\mathcal{L}$, namely $\mathcal{L}_n(\{0\})\leq 1- \kappa$ for all $n$ with some $\kappa>0$ (thus ensuring that at least one leaf is added with probability uniformly bounded away from $0$). The \name{} has been further studied  \cite{englander2021tree}, dropping the uniform ellipticity assumption. 

We refer to this model as $\mathcal{L}$-\name\  to emphasize the  dependence on the sequence $\mathcal{L}:=\{\mathcal{L}_n\}_{n\ge 1}$, which accounts for different probabilities of adding new vertices to the tree along the evolution. We denote by $\Pd_{T_0, x_0; \mathcal{L}}$ the law of $\{(T_n, X_n)\}_{n\geq 0}$ when $(T_0, X_0)=(T_0,x_0)$  and by $\Ed_{x_0, T_0;\mathcal{L}}$ the corresponding expectation.
Given a sequence of laws $\mathcal{L}:=\{\mathcal{L}_n\}_{n\ge 1}$ and a natural number $m$, we will write $\mathcal{L}^{(m)}$ for the shifted sequence of laws, that is, $\mathcal{L}^{(m)}: = \{\mathcal{L}_{m+n}\}_{n\ge 0}$.

\subsection{The Barab\'asi-Albert model}\label{ss:BA}

The Preferential Attachment Random Graph Model proposed by A-L. Barabási and R. Albert in \cite{barabasi1999emergence}, explains the formation of scale-free networks, characterized by a few highly connected nodes. The model is based on the principle that new nodes are more likely to connect to already well-connected nodes, leading to ``hubs''. Its rigorous mathematical theory began with the seminal paper \cite{bollobas2001degree} and  was further developed by B. Bollob\'as and O. M. Riordan in a series of papers; see, e.g., \cite{BollobasRiordan2003}. 

The classical version of the model evolves as follows. Start with an initial (finite) graph $G_0$ (not necessarily a tree).
Then, obtain $G_{t+1}$ from $G_t$ in the following way: add a vertex $v_{t+1}$ to $G_t$ and connect it to a vertex $u \in G_t$ selected independently of the past and with probability proportional to its degree. More formally, conditioned on $G_t$, the probability of connecting $u$ to the new vertex $v_{t+1}$ is
\begin{equation}\label{eq:pa}
    P\left( v_{t+1} \text{ connects to }u \; \middle |\; G_t\right) = \frac{\text{degree of } u \text{ in }G_t}{\text{sum of all degrees in }G_t}.
\end{equation}
The above equation is the so-called {\it preferential attachment rule}, in which new vertices are more likely to connect themselves to ``popular'' vertices (i.e. vertices with higher degree). Throughout this paper, we will refer to the above defined model as {\it BA-model}.   The BA-model is sometimes considered starting from a single vertex with a self-loop (henceforth denoted by $G_0^{\rm loop}$) and  sometimes starting from  a single edge between two vertices ($G_0^{\rm edge}$). 

When $G_0$ is a tree or a tree with a unique vertex with a self-loop, we shall refer to the model as BA-tree (despite the slight abuse of terminology).

There is a vast literature about this rule and its many variations. We refer the reader to \cite{remcoVol1,remcoVol2} for a comprehensive discussion of preferential attachment random graphs.

 \subsection{Applications}\label{ss:app}
Before presenting our main results, we discuss some interesting applications of the TBRW.

\subsubsection{Chain Walk Catalyst} 
Polyethylene is one of the most widely produced and used types of plastic, with many different applications. Depending on how it is produced, it can have different structures, such as Low Density Polyethylene, Linear Low Density Polyethylene, and High Density Polyethylene. These types are classified based on their density and branching patterns. However, creating a specific and precise polymer structure is more difficult.

In \cite{guan1999chain,johnson1995new} the authors introduced a special catalyst to make polyethylene and similar materials with various structures depending mostly on the pressure, used during synthesis. This catalyst has the unique ability to attach to the polymer and randomly walk along its structure, creating complex branching patterns, a process known as {\it ``Chain Walking''} (CW) catalysis. Over time, it has been shown through experiments and simulations that the CW mechanism can produce linear-like chains at high pressure, while highly branched tree-like structures at low pressure.

The \name{} model provides a mathematical framework that can be seen as a model for this type of catalyst. Since in the \name{} the walker moves step-by-step, either creating new branches or extending existing ones, this process mirrors the way the CW catalyst moves along the polymer chain. Also, the TBRW model allows for the generation of a variety of tree-like structures, similar to the branching patterns observed in polymers synthesized using CW catalysis. To model different 
 synthesis conditions, such as pressure, the TBRW model can generate different tree structures by choosing different sequences of laws $\mathcal{L}$. High pressure conditions in CW catalysis lead to more linear chains, akin to a TBRW with 
$\mathcal{L}_n = \mathsf{Ber}(p)$ for all $n$, whereas low pressure conditions lead to highly branched structures, similar to a TBRW with $\mathcal{L}_n = \mathsf{Ber}(p_n)$ where $p_n$ goes to zero 
fast enough. 

Our paper is a follow up to  \cite{englander2021tree}, where the focus was more on the random walk than the graph structure. Since \cite{englander2021tree}  has recently been cited in \cite{DS22}, a paper published in the Journal of Chemical Physics, we hope that the results in this current paper are also relevant in chemical physics.

\subsubsection{On the ubiquity of linear    preferential attachment and an approach to it via TBRW}

The BA model provides a mechanism for generating networks with power-law degree distributions, a feature commonly observed in numerous real-world networks such as the Internet, biological nets, transportation networks, power grids and many others~\cite{10.1145/316194.316229,10.1093/acprof:oso/9780198515906.001.0001,BRODER2000309}. The preferential attachment rule behind the BA model has helped us better understand the structure of the World Wide Web and the emergence of highly popular web pages \cite{barabasi1999emergence, RevModPhys.74.47}.

However, it is not entirely clear why the assumption of preferential attachment (PA) should hold---namely, why new vertices are more likely to connect to pre-existing vertices with higher degrees, especially because this means that incoming vertices must know the degrees of all the other vertices in the network.  It is even less clear why preferential attachment should be {\it linear}, i.e., proportional to the degrees, maybe with an additive constant, but not proportional to a power of the degrees, say. Non-linear PA rules do not yield power-law degree sequences \cite{krapivsky2000connectivity,oliveira2005connectivity,rudas2007random}, hence the ubiquity of power-laws needs a more robust explanation. 

The idea that random graphs with a power-law degree distribution should be built locally by random walks was first proposed in \cite{vazquez2003growing,saramaki2004scale,evans2005scale}, with numerous possible real-world applications listed there, from brain development to technological innovation. However, in all these models the random walk runs for a few steps before new vertices are added and in each growth phase the starting position is resampled from the uniform or other distribution. Rigorous mathematical treatments of such models have been undertaken in \cite{CJ13, banerjee2024co, berry2024random}. Such growth mechanisms can be achieved locally if every existing vertex has an independent clock for the creation of a new vertex (with a rate possibly depending on the local neighborhood), but, if we have only one or a bounded number of agents who can grow the graph (as in the CW catalysis example of the previous subsection), then choosing the growth place from the uniform or any other distribution on the entire graph is not a local procedure. Graphs generated by random walks without restarts, as in the present paper, were first proposed in \cite{amorim2016growing,Iacobelli_Figueiredo_Neglia_2019} and \cite{ross2019random}. 

\subsubsection{Physical networks.} In the recent concept of a physical network \cite{dehmamy2018structural,posfai2024impact}, not only is the graph embedded in space, but the vertices and edges are non-overlapping physical objects. In \cite{pete2023network}, a network-of-networks model for physical networks was introduced: a randomly evolving subgraph of an ambient graph (say, a random walk trajectory in a box of $\mathbb{Z}^d$) grows until it hits the already existing structure, building a tree from how these spatially extended nodes touch each other. This is very relevant for {\it modeling brain networks}, for instance, where large neurons connect to each other through point-like synapses, or {\it river networks}, and so on. It was observed non-rigorously in \cite{pete2023network} that, when the ambient dimension is more than twice the fractal dimension of the randomly growing pieces, then the place of attachment for the new piece is well-mixed in space, hence depends mainly on the sizes of the existing pieces, and a linear preferential attachment mechanism arises. This is quite similar to the mechanism in the present paper, except that the random motion in \cite{pete2023network} takes place in the ambient space, while in the present paper it takes place inside the existing structure; a situation analogous to Diffusion Limited Aggregation (DLA) \cite{witten1983diffusion} versus Internal DLA \cite{lawler1992internal}. It would be interesting to study what kind of physical networks could arise from internal growth.

\subsection{Main results}\label{ss:res}

In order to state our main result and discuss its consequences, we need to introduce some definitions.

Let $\mathsf{G_{fin}}$ denote the set of finite graphs,  $G_0 \in \mathsf{G_{fin}}$ and  $P_{G_0}$ the law of the BA-model, conditioned on the initial graph being $G_0$. Given a 
family $\mathcal{G}$ of sequences of graphs in $\mathsf{G_{fin}}$ we say that  a random graph sequence $\{G_t\}_{t \in \mathbb{N}}$ satisfies $P_{G_0}$-almost-surely the property $\mathcal{G}$,  if $P_{G_0}(\{G_t\}_{t \in \mathbb{N}} \in \mathcal{G})=1$. 
\begin{definition}[Asymptotic Graph Property] We say that a family $\mathcal{G}$ of sequences of graphs is an {\rm asymptotic graph property} if satisfying (belonging to) $\mathcal{G}$ does not depend on a finite number of coordinates. More formally, it means that for any $s\in\mathbb{N}$, we have  
$$
\{H_t\}_{t\in\mathbb{N}} \in \mathcal{G} \iff \{H_{t+s}\}_{t\in\mathbb{N}} \in \mathcal{G}.
$$
\end{definition}
\begin{remark}
Note that if $\mathcal{G}$ is an asymptotic graph property and $P_{G_0^{\rm loop}}(\{G_t\}_{t \in \mathbb{N}} \in \mathcal{G})=1$ then, by the Markov property,  $P_{G}(\{G_t\}_{t \in \mathbb{N}} \in \mathcal{G})=1$, for all $G$ finite tree with a unique vertex having a self-loop.       
\end{remark}

We give three examples to illustrate the concepts of asymptotic and almost sure graph properties.

\begin{example}[Linear growth] Consider $\mathcal{G}$ as the following property
    $$
    \mathcal{G} = \left\lbrace\{G_t\}_t \, : \, \lim_{t\to \infty}\frac{|G_t|}{t} = 1\right \rbrace,
    $$
    where $|G_t|$ denotes the cardinality of the vertex set of $G_t$. Then, $\mathcal{G}$ contains those sequence of graphs whose vertex set grows linearly in time, with speed 1. Then, $\mathcal{G}$ is an asymptotic graph property. For the BA-model, it is  a $P_{G_0}$-almost sure property.
\end{example}

\begin{example}[The height of the BA-tree] For a fixed positive real number $h>0$, consider the property
$$
\mathcal{G}_h = \left\lbrace\{G_t\}_t \, : \, \lim_{t\to\infty}\frac{\text{\rm height of }G_t}{\log |G_t|} = h\right \rbrace.
$$
Then, $\mathcal{G}_h$ is an asymptotic property for any $h$. On the other hand, \cite[Theorem 1]{pittel1994note} guarantees (see also \cite[Section 8.2]{remcoVol2}) that there exists a positive constant $h_0$ for which, $P_{G_0^{\rm edge}}$-almost surely,  the BA-tree satisfies $\mathcal{G}_{h_0}$. 
\end{example}

\begin{example}[The maximum degree of the BA-tree] For any fixed real number $h\ge 0$, consider the property
$$
\mathcal{G}_h = \left\lbrace\{G_t\}_t \, : \, \lim_{t\to\infty}\frac{\mathsf{max.deg}(G_t)}{|G_t|^{1/2}} > h\right \rbrace.
$$
Then, $\mathcal{G}_h$ is an asymptotic property for any $h$. It is shown in \cite[Theorem 3.1]{mori2005maximum} that it has a $P_{G_0^{\rm edge}}$ non-trivial probability for the BA-tree for any $h>0$, hence it is not a $P_{G_0^{\rm edge}}$-almost sure  graph property.
(The notation in the denominator refers to the maximal degree among the vertices of $G_t$.)
\end{example}

Our main results give us conditions under which we can transfer asymptotic graph properties which hold almost surely for the BA-tree to the graph sequence generated by the \name; we call this transfer mechanism the {\it Transfer Principle}. In this paper, we introduce a mixing condition~\mixc. 
Our main result states that, under \mixc, the Transfer Principle holds.

Before proceeding, we need to introduce a sequence of stopping times $\{\tau_k\}_{k \in \mathbb{N}}$ corresponding to the times when the \name{} process actually adds at least one new vertex, the {\it growth times}. The sequence is defined  as follows: $\tau_0 \equiv 0$ and for $k\geq 1$ 
\begin{equation}\label{def:tauk}
    \tau_k := \inf \{n>\tau_{k-1} \, : \, Z_n \ge 1\},
\end{equation}
if $\tau_{k-1} <\infty$, and $\tau_k = \infty$ otherwise. Observe that the sequence $\{T_{\tau_k}\}_{k \in \mathbb{N}}$ is the subsequence of $\{T_n\}_{n\in \mathbb{N}}$ that carries all modifications made by the walker. More formally, we have that $T_n = T_{\tau_k}$ for all~$n \in [\tau_k,\tau_{k+1})$. Throughout this paper we will consider only sequences of probability laws $\{\mathcal{L}_n\}_{n\ge 1}$ under which $\tau_k < \infty$ holds for all $k$ almost surely. This is because we are not interested in scenarios in which the graph sequence will eventually be constant (stop growing) with positive probability. On the other hand, our mixing condition \mixc{}  below will ensure, roughly speaking, that the growth times are rare enough so that the random walk on each of the current trees (after finitely many exceptions) mixes before a new vertex is added.

 Let $\Delta\tau_k:=\tau_{k+1}-\tau_k$. Furthermore, for each $k\in \mathbb{N}$, let $\eta_{k}$ be an optimal strong stationary time for the TBRW started at~$(T_{\tau_{k-1}},X_{\tau_{k-1}})$ and with the shifted sequence of laws $\mathcal{L}^{(\tau_{k-1})}$, as specified later in Section~\ref{ss:ssttbrw}. (We do not spell it out here to avoid too many definitions in this introduction.)

\begin{definition}[Condition \mixc{}]\label{def:TP} We say that  $\mathcal{L}$-\name{} satisfies \mixc{} if the following conditions hold:
\begin{enumerate}
    \item[(1)] For all $n$, $\mathcal{L}_n=\mathrm{Ber}(p_n)$, with $p_n\in (0,1)$.
    \item[(2)]  For any finite tree initial condition $(T,x)$, there exists some sequence of optimal strong stationary times $\{\eta_k\}_k$ such that
        $$
            P_{T,x;\mathcal{L}^{}}\left(\eta_k \geq \Delta\tau_k \text{ i.o.} \right) = 0.
        $$
\end{enumerate}
\end{definition}
In words, item (1) means that at each step the walker can add at most one vertex, and has a positive probability both for adding one and for adding none. This is necessary for two reasons: firstly, when $\mathcal{L}_n(\{1\})=1$, it was proved in \cite{10.1214/20-EJP574} that the limit random tree generated by the \name{} is one-ended, which is not true for the BA model. Secondly, we need to rule out the pathological case $\mathcal{L}_n(\{1\})=0$ when the walker is just a SSRW over the initial condition.  Item (2) of \mixc{} is a mixing condition. It guarantees that after some random time, all subsequent new vertices are added after the walker has mixed on the current tree.

\begin{theorem}[Transfer Principle under \mixc]\label{thm:generaltransfer} Let $\mathcal{G}$ be an asymptotic graph property. Consider a \name\ satisfying condition \mixc\ of Definition~\ref{def:TP}. Then, 
\[
P_{G_0^{\rm loop}}\left(\{G_t\}_{t \in \mathbb{N}} \in \mathcal{G}\right)=1 \implies \mathbb{P}_{T_0,x_0;\mathcal{L}^{}}\left( \{T_{\tau_k}\}_{k \in \mathbb{N}}\in \mathcal{G}\right)=1, \quad \text{ for all $(T_0,x_0)$}.
\]
\end{theorem}

The next result provides a large class of sequences of distributions $\mathcal{L} = \{\mathcal{L}_n\}_{n\geq 1}$ for which the $\mathcal{L}$-\name{} satisfies condition \mixc.

\begin{theorem}[A family of \name\ satisfying \mixc]\label{thm:transfer} The $\mathcal{L}$-\name{} with $\mathcal{L}_n = {\rm Ber}(p_n)$ and $p_n \asymp n^{-\gamma}$ (meaning that there exist uniform constants $0<\g,\G<\infty$ with $\g n^{\gamma} < p_n < \G n^{-\gamma}$ for all $n$) satisfies condition~\mixc\ for all~$\gamma \in (2/3,1]$.
\end{theorem}

The following results are consequences of our Transfer Principle and give us several structural results about the random graph sequence $\{T_n\}_{n \in \mathbb{N}}$ generated by the \name.

\begin{corollary}[Power law for degree distribution]\label{cor:degdist} Consider a $\mathcal{L}$-\name\ satisfying condition \mixc. Then, for any finite initial condition $(T,x)$ 
and degree $d \geq 1$, we have that
$$
 \lim_{n\to \infty}\frac{\# \text{ \rm of vertices having degree exactly $d$ in }T_n}{|T_n|} =  \frac{4}{d(d+1)(d+2)},
$$
$\Pd_{T,x; \mathcal{L}^{}}$-almost surely.
\end{corollary}

This corollary was proved in \cite[Theorem 1.5]{englander2021tree} for the particular sequence of laws~$\mathcal{L}_n = \mathsf{Ber}(n^{-\gamma})$ with $\gamma \in (2/3,1]$, using stochastic approximation methods. In our case we not only have a more general statement, but also give a proof that is simpler and conceptually more appealing.

\begin{corollary}[Height of $T_n$]\label{cor:height} Let $c$ denote the unique solution of the equation $ce^{c+1} =1$. Then, for any $\mathcal{L}$-\name\ satisfying condition \mixc, and
any finite initial condition $(T,x)$, 
 $$\lim_{n\to \infty}\frac{\text{\rm height of }\, T_n}{\log |T_n|} =  \frac{1}{2c},
$$
$\Pd_{T,x; \mathcal{L}^{}}$-almost surely.
\end{corollary}

\begin{corollary}[Maximum degree]\label{cor:maxdeg} Consider an $\mathcal{L}$-\name\ satisfying condition \mixc. Then, for any finite initial condition $(T,x)$,  
there exists a strictly positive  random variable $\zeta$ such that
$$
\lim_{n\to \infty}\frac{\mathsf{max.deg}(T_n)}{\sqrt{|T_n|}}  = \zeta, 
$$
$\Pd_{T,x; \mathcal{L}^{}}$-almost surely. 
\end{corollary}

We finish the list of corollaries with  a scaling limit result which we will derive from a result in \cite{mori2006surprising}, using the Transfer Principle.

We say that a  vertex, different from the root, is the {\it child} of another vertex  if it is the direct descendant of the latter  in the tree. The {\it zero level} of the tree  consists of the root only. Children of the root
 form the {\it first level}, children of vertices at first level form the {\it second level}, and so forth. 

 Now focus on the ``lower'' levels, i.e., fixed (not depending on $n$) levels, different from the zero level. The following result gives the size of the $k$th level for a fixed $k\geq 1$, and characterizes the degree distribution of vertices belonging to the $k$th level. It reveals the surprising fact, that the asymptotic degree distribution of vertices on lower levels is {\it different} from the one when the whole tree is considered, and of the largest levels. Although it still follows a power law, the exponent is not $3$ as in Corollary \ref{cor:degdist}, but $2$, for each level $k\geq 1$.

Let
$X[n,k]$ denote the size of level $k$ after step $n$ and let
$X[n, k,d]$ denote the number of vertices at level $k$ with degree $d$, after step $n$. 
\begin{corollary}[Degrees of $k$th level vertices.]\label{cor: kth.level}
Consider an $\mathcal{L}$-\name\ satisfying condition \mixc. Then, for any finite initial condition   $(T,x)$  
and for every $k\ge 1,\ d\ge 1$, 
$$\lim_{n\to\infty}\frac{X[n,k,d]}{X[n,k]}=\frac{1}{d(d+1)},$$ $\Pd_{T,x; \mathcal{L}^{}}$-almost surely. (cf.  Corollary \ref{cor:degdist}).
\end{corollary}

We conclude this subsection with five (groups of) open problems.

\begin{problem}[The walker's height process]
For $\cL_n=\mathsf{Ber}(p_n),$ let $H_n$ denote the distance of the walker from the root (the ``height'' of the walker) at time $n$. By Corollary~\ref{cor:height}, if Condition~\mixc\ is satisfied, 
then the diameter (hence the height of the tree) scales with $\log |T_n|$ a.s.

In order to gain some insight about the process $H$, it will  be useful to recall two results on BA trees due to Z.~Katona. Define $m_n=\sigma_n^2=(1/2)\log n,\, n\ge 1$. Let  $\varphi$ denote the standard normal density:
$$\varphi(r):=(\sqrt{2\pi})^{-1}e^{-r^2/2}.$$
The first result is \cite[Theorem 1]{Katona05}, giving a Local Central Limit Theorem for the level distribution of the population.
Namely, denoting by $X[n,k]$, as before, the population size at level $k\ge 0$, one has almost surely (in the randomness of the BA process) that
\begin{align}\label{Katona.localCLT}
\lim_{n\to\infty}\sup_{r\in K}\left |\sigma_n\frac1n X\big[n,\, \lfloor m_n +r\sigma_n\rfloor \big]-  \varphi(r)\right|=0,
\end{align}
for any compact interval $K\subset \mathbb R$. Loosely speaking, the empirical distribution of the heights of the vertices is close to being normal with mean $m_n$ and variance $\sigma_n^2$, when $n\gg 1$.

Katona's second result \cite[Theorem 1]{Katona06} says that the degree distribution on any level $\lfloor m_n +r\sigma_n\rfloor$ agrees asymptotically with the (exponent 3 power law) degree distribution of the entire tree. 
Putting this together with~\eqref{Katona.localCLT} shows that the bulk of the stationary distribution of a  simple random walk on the BA tree is also at levels that are in the $O(\sigma_n)$-neighborhood of level $m_n$. More precisely, consider $$\pi[n,k]:= \frac{1}{2n} \sum_{d\ge 1} d \, X[n,k,d],$$
the total stationary measure of the height $k$ vertices in a BA tree with $n$ vertices (with a loop at the root). Then the above results \emph{suggest} the following analog of \eqref{Katona.localCLT}: almost surely,
\begin{align}\label{Katona.piCLT}
\lim_{n\to\infty}\sup_{r\in K}\left |\sigma_n\pi\big[n, \lfloor m_n +r\sigma_n\rfloor \big] \,   -  \varphi(r)\right|=0,
\end{align}
for any compact interval $K\subset \mathbb R$. However, the results of~\cite{Katona06} are not quite strong enough to deduce this; for instance, a few vertices of large degrees would not ruin the asymptotic degree distribution, but could modify the stationary distribution.

Returning now to our model, by virtue of the coupling, the tree at time $n$ is close to a BA tree at time $\widehat n:=|T_n|.$ (Note that $|T_0|=1$.)
Even though $\widehat n$ is random,  it can be determined asymptotically a.s., using the (extended) Law of Iterated Logarithm. 
For example, when $p_n \asymp n^{-\gamma}$ and $\gamma\in (0,1]$, we have that $\widehat n$ is approximately
$n^{1-\gamma}$. 
Under Condition~\mixc, i.e., for $\gamma \in (2/3,1]$, the TBRW spends most of its time in stationarity on the current tree; in fact, at any large deterministic time $n$, the probability of being so close to a growth time that we have not had enough time yet to get very close to stationarity on the new tree is polynomially small in $n$. This is negligible compared to the stationary distribution $\asymp 1/\sqrt{\log{\hat{n}}} \asymp 1/\sqrt{\log{n}}$ of any given height in the scaling window given in~(\ref{Katona.piCLT}). Thus, the asymptotic results above \emph{suggest} that the one-dimensional distribution of the height process satisfies the local CLT
\begin{equation}\label{Hlim}
\lim_{n\to\infty}\sup_{r\in K}\left |\sigma_{\hat{n}}\,\Pd_{T_0, x_0; \mathcal{L}} \big( H_n = \lfloor m_{\hat{n}}+r\sigma_{\hat{n}}\rfloor \big) \,  - \varphi(r)\right|=0\,,
\end{equation}
again, for every compact interval $K\subset\R$. (A minor complication when using Katona's results is that he does not seem to have a loop at the root; this is unlikely to drastically change anything, but warrants an additional step in a proof.)

A slightly stronger statement concerns the averaged time behavior. 
Consider the local time $\ell$ of the process $H$ defined for level $k$ at time $n$ as
$$\ell(n,k):=\sum_{i=0}^n \mathbb 1 _{\{H_i=k\}}.$$
Since the walker mixes rapidly between growth times, it is natural to guess that $\ell$ (normalized by $n$) is close to the stationary distribution. More precisely, if the LCLT~\eqref{Katona.piCLT} for the stationary distribution holds, then it should also be the case that 
\begin{align}\label{HellCLT}
 \lim_{n\to\infty} \sup_{r \in K} \left|\frac{\sigma_{\hat{n}}}{n} \ell\big(n,\lfloor m_{\hat{n}}+r\sigma_{\hat{n}}\rfloor \big) - \varphi(r)\right| = 0
\end{align}
for any compact $K\subset\R$, in probability, or even almost surely. This could follow from the quantitative Markov Chain Ergodic Theorem \cite[Theorem 12.21]{levin2017markov} applied to the random walk on the ``frozen'' graphs between growth times.

The next natural question would be the full (functional) scaling limit of the processes $(H_{nt})_{t\ge 0}$  where, as usual, the process at integer times  is extended to all real times by linear interpolation.  As a direct continuation of our conjectures~\eqref{Hlim} and~\eqref{HellCLT}, we first look at the time scale $\sigma^2_{\hat{n}} \asymp \log n$, much shorter than the total time $n$. As a natural comparison, one should also look at possible scaling limits of the height process $(H^{\mathrm{BA},n}_{t})_{t\in [0,\infty)}$  (again, the process at integer times is extended to all real times by linear interpolation) of a simple random walk on a BA tree with fixed $n$ vertices, started from the stationary distribution (corresponding to the fact that the \name{} spends most of its time in stationarity). In light of~\eqref{Katona.localCLT} and~\eqref{Katona.piCLT}, one may guess that \emph{if} the centered and scaled processes
\begin{align}\label{BAtlogn}
\mathcal{H}_n^{\mathrm{BA}}(t):=\frac{H^{\mathrm{BA},n}_{t\sigma^2_n} - m_{n}}{\sigma_{n}}, \quad t\in [0,\infty)
\end{align}
have a Markovian scaling limit as $n\to\infty$, then the limit should be a stationary Gaussian process, such as an Ornstein-Uhlenbeck process. (We are not aware of such results in the literature.) 
Then, for the \name{}, since at a typical moment $n$ the process is close to stationarity on its current tree, and far from adding a new leaf, it is natural to guess that the processes
\begin{align}\label{Htlogn}
\mathcal{H}_{n}(t):=\frac{H_{n+t\sigma^2_{\hat{n}}} - m_{\hat{n}}}{\sigma_{\hat{n}}}, \quad t\in [0,\infty)
\end{align}
have the same limit as~\eqref{BAtlogn}.

As far as the \emph{difference} between $(H_{nt})_{t\in [0,\infty)}$ and $(H^{\mathrm{BA},\hat n}_{\hat{n}t})_{t\in [0,\infty)}$ is concerned, that could become visible only on time-scales much larger than the scale $\sigma_n^2 \asymp \sigma_{\hat n}^2 \asymp \log n$ that we have been discussing so far. Namely, since \name\  grows a new leaf only at every $\approx n^\gamma$th step (at a crude, undefined level of precision), the time scale should not be smaller than this if we want to see the growth. However, we expect the hitting and mixing times on a BA tree with $\hat n \approx n^{1-\gamma}$ vertices to be $n^{1-\gamma+o(1)}$, which is of much smaller order. That is, on the time scale $n^\gamma$, the walker is ``all the time everywhere on the tree simultaneously,'' hence there seems to be no interesting scaling limit for the height of the TBRW on this time-scale.
\end{problem}

\begin{problem}[The regime $1/2 < \gamma\leq 2/3$ for $p_n \asymp n^{-\gamma}$]
Here it still holds that an optimal stationary time is typically reached earlier than the next growth time, hence the growth process is typically similar to the BA tree process; however, the growth a.s.~happens infinitely many times much earlier than any optimal stationary time; in fact, w.h.p.~there are at least $n^{\eps}$ growth times that are at most $n^{-\eps}$ fraction of the volume of the current tree, for some $\eps=\eps(\gamma)>0$. Based on this, we are expecting that the process is actually distinguishable from the BA process: their total variation distance, up to time $n$, goes to $1$ as $n\to\infty$.

The intuition is as follows. In the BA tree process, it is easy to prove that within the first $n$ steps w.h.p.~there are only $n^{o(1)}$ instances where the smallest subtree containing two consecutive leaf additions has volume $n^{o(1)}$. On the other hand, if the TBRW process could be coupled with a uniformly positive probability to agree with the BA process from some time on, then one should be able to prove, using known structural properties of BA trees, that, among the growth times that are at most $n^{-\eps}$ fraction of the current volume, there are $n^{\eps'}$ many, for some $\eps'>0$, in which the walker did not manage to get far enough to escape a small sub-tree containing the previous leaf addition --- a contradiction.

We will not try to make this argument rigorous, since this is a negative result of limited significance: we  expect that the interesting large scale statistics of the TBRW tree (asymptotic degree distribution, diameter, etc.) still  agree with those of the BA tree.
\end{problem}

\begin{figure}[htbp]
\begin{tabular}{p{4.3 cm}p{1.3cm}p{4.3 cm}p{1.2cm}}
\includegraphics[width=4.5 cm]{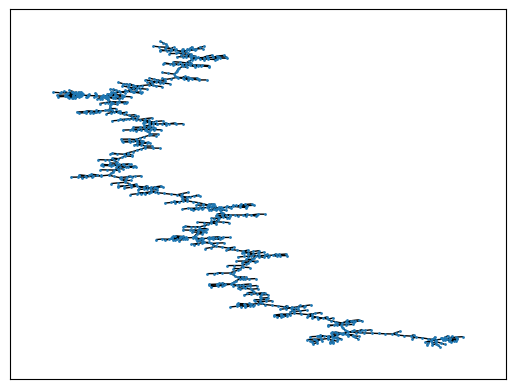} &
\vskip - 1cm $\gamma=0$ &
\includegraphics[width=4.5 cm]{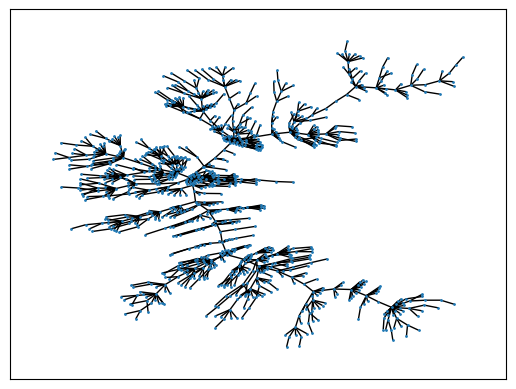} &
\vskip - 1cm  $\gamma=0.3$\\
\includegraphics[width=4.5 cm]{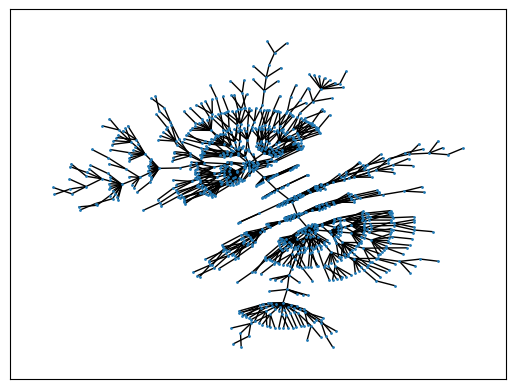} &
\vskip - 1cm  $\gamma=0.55$ &
\includegraphics[width=4.5 cm]{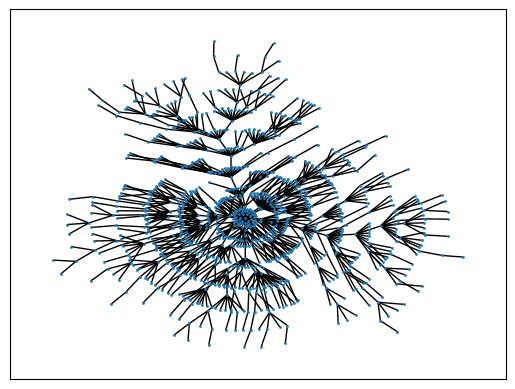} &
\vskip - 1cm  $\gamma=0.9$ 
\end{tabular}
\caption{Simulation of TBRW with 1000 vertices, with different values of $\gamma$. (Thanks to Á. Kúsz for the pictures.)}
\end{figure}

\begin{problem}[The regime $0 < \gamma < 1/2$ for $p_n \asymp n^{-\gamma}$]
Here the times between leaf additions still go to infinity, but are much smaller than the mixing time of the current tree. Hence new leaves are typically added far from the previous addition 
and also far
from stationarity. Thus we expect that the diameter of the tree grows polynomially but sublinearly in its volume, like $n^{\delta(\gamma)+o(1)}$ for some $0<\delta(\gamma)<1-\gamma$, and the walk is transient: any finite subset of the vertices is visited only finitely many times. We do not have a guess for the degree distribution or for the shape of the local limit of the tree as viewed from the walker (shown for $\gamma=0$ to be one-ended, with exponentially decaying degree distribution, in \cite{10.1214/20-EJP574}).
\end{problem}

For $\gamma=1/2$, we do not even have a guess, regarding recurrence vs.~transience.

 \begin{problem}[More general rates $(p_n)$] In this paper we investigate the \name{} with $\mathcal{L}_n = {\rm Ber}(p_n)$ and $p_n \asymp n^{-\gamma}$. It was suggested to us by one of the referees to consider more general sequences $p_n \to 0$. Is there a general summability type of condition to make sure that condition~(M) of Definition~\ref{def:TP} holds?

One strategy could be trying to show that a graph observable like the diameter is ``non-increasing'' in $(p_n)$, meaning that if $(p_n)$ and $(p'_n)$ are the parameters of two Bernoulli sequences, with the former decreasing faster, then $\mathsf{diam}(T_n) \le \mathsf{diam}(T'_n)$ should hold. This way, the walker of a TBRW with a faster rate would mix faster, consequently satisfying $(M)$ more easily.  However, showing this kind of monotonicity of a graph observable that depends on the geometry of $T_n$ and on how the walker explores its environment is hard in this model, due to the lack of monotone coupling between models with different parameters. 

What should be possible is to rewrite the \emph{proofs} to include the more general condition that $p_n$ is a decreasing sequence satisfying $c n^{-1} < p_n < Cn^{-\gamma}$ for some $\gamma>2/3$ and $0<c,C<\infty$ fixed. We have chosen not to complicate the writing further with this generality. 

An exact summability criterion for $(p_n)$ seems currently to be out of reach, due to our use of the suboptimal \emph{a priori} diameter bound Theorem~\ref{thm:diam} and the mixing time bound Corollary~\ref{cor:mixingbound}.
\end{problem}

\begin{problem}[More general random walks] 
One could consider versions of the \name{} where the walker is moving not according to simple random walk, but some other nearest-neighbor Markov chain, such as a bias towards the root (as in \cite{banerjee2024co}, for instance), or the \emph{degree-biased random walk}, popular in network science \cite{gomez2008entropy,sinatra2011maximal,bonaventura2014characteristic}. On a tree, the latter is easiest to describe as the reversible chain with electric conductances $c(x,y):=\deg(x)^\alpha\deg(y)^\alpha$ on the edges $(x,y)$, for some $\alpha\in [0,\infty)$. For any such random walk, our strategy going through strong stationary times and the transfer principle of Theorem~\ref{thm:generaltransfer} could be attempted in order to find how the trees built by this walk look like asymptotically.
\end{problem}

\subsection{Proof ideas and technical difficulties}\label{ss:ideas}

In a nutshell, the proof of Theorem \ref{thm:generaltransfer} follows from a coupling with the BA-model starting from a random initial condition. In order to construct such a coupling, we use the framework of optimal strong stationary times.

Regarding Theorem \ref{thm:transfer}, the main idea behind its proof is to guarantee that, when we set $\mathcal{L}_n = {\rm Ber}(p_n)$, if the sequence $(p_n)_{n\geq 1}$ goes fast enough to zero, the walker mixes before adding new vertices (this is item (2) of condition \mixc). This guarantees that new vertices are added according to the stationary measure, which allows us to couple the \name\ with an instance of the BA model.

Although the core ideas behind our proofs are intuitive, their implementation is far from straightforward. The two main issues to overcome are:
\begin{enumerate}
    \item Good enough bounds for the mixing time depend on structural properties of the graph, such as its diameter.
    \item Mixing on the current tree before adding new vertices ensures that new vertices are added according to a distribution that is close to the stationary distribution of the tree; however, since we are not exactly in the stationary distribution, this is not enough to couple the model with the BA model. 
\end{enumerate}
 
As far as issue (1) is concerned, it is key to guarantee a good enough upper bound for the mixing time, since a small mixing time increases the chances of mixing before adding a new vertex. The easiest upper bound for mixing time to obtain is the number of vertices squared, see Corollary \ref{cor:mixingbound}. Since the order of the graph at time $n$ is a sum of independent random variables, we have good control over this quantity. However,  depending on $(p_n)_{n\geq 1}$, we may observe infinitely often the walker adding a new vertex before taking this many steps, and our purpose is to find a better bound on the mixing time and conditions that rule  out adding vertices before mixing. This leads us to a stronger general bound for mixing time on trees, which is basically the number of vertices times the diameter of the tree; see Corollary~\ref{cor:mixingbound} below. This bound represents a substantial improvement as long as we have good enough control over the diameter's growth. Here is where we face one of the main challenges in the analysis of the \name: {\it we cannot decouple the walker's trajectory from the environment's evolution}. Thus, in order to obtain an upper bound on the diameter, one needs to control the walker's trajectory. On the other hand, the walker's trajectory depends on the underlying graph. 

We dedicate Section \ref{sec:diameter} to proving a sub-polynomial upper bound on the diameter (Theorem~\ref{thm:diam}). This bound is good enough to allow us to couple the \name\ with the BA model and consequently obtain the right order of the diameter, which is logarithmic in the number of vertices.

Regarding difficulty (2), it is important to highlight the fact that the walker having mixed on a graph does not imply that it is distributed as the stationary distribution. However, to couple  the \name\ with the BA model, we need that the new leafs added by the process are added at a position chosen according to the stationary distribution. This leads us to the notion of stationary times. Dealing with stationary times comes with extra challenges as well. Moreover,  we also need to ensure that the random times at which new vertices are added after the stationary times remain stationary times. We refer the reader to Section~\ref{sec:prep}, where we introduce and discuss all the key concepts needed throughout the paper.

\subsection{Organization of the paper}

This paper is organized as follows:
\begin{itemize}
    \item In Section \ref{sec:prep}, as a preparation, we introduce the concepts of mixing and (strong) stationary times together with some results about these topics. We also discuss the concept of stationary times in the context of the \name, since our trees are changing over time. 
    \item In Section \ref{sec:generaltransfer}, we prove Theorem~\ref{thm:generaltransfer}, which guarantees the Transfer Principle under condition \mixc.
    \item In Section \ref{sec:transfer}, we prove Theorem \ref{thm:transfer}, assuming the subpolynomial bound Theorem~\ref{thm:diam} on the diameter.
    \item In Section \ref{sec:diameter}, we prove our subpolynomial upper bound on the diameter.
\end{itemize}


%
%

\section{Preliminaries}\label{sec:prep}
In this section we will set the stage by introducing all the key concepts needed in the proof of our main theorem. 

\subsection{Mixing time upper bound}

Let us begin by recalling same basic facts about the spectral representation of Markov chains (see \cite[Section 12.1]{levin2017markov} for a more detailed account). 
Let $\{Y_t\}_t$ be an irreducible Markov chain on a finite state space $S$ with transition matrix $P$.  We assume that $P$ is reversible with respect to the probability distribution $\pi$, i.e., $\pi(x)P(x,y)=\pi(y)P(y,x)$, for all $x,y \in S$ (thus, $\pi$ is the stationary distribution). We point out here that in this paper we only apply the general theory to  Markov chains  that are reversible.

Since $P$ is irreducible and reversible, 
it has a complete spectrum of real eigenvalues with magnitude at most $1$. These eigenvalues are denoted by
\[
1=\lambda_1 >\lambda_2\geq \cdots \lambda_{|S|}\geq -1,
\]
where $\lambda_2<1$ by irreducibility. The chain is aperiodic if and only if  $\lambda_{|S|}>-1$.  

The {\it absolute spectral gap}, denoted by $\gamma_*$, is defined as 
$\gamma_*:=1- \max\{|\lambda_i|: i\geq 2\}$. The {\it relaxation time} of a reversible Markov chain with absolute spectral gap $\gamma_*$ is defined as 
\[
t_{\rm rel}:=\frac{1}{\gamma_*}.
\]
For $(1/2)$-lazy Markov chains (where we take a step of the chain only with probability $1/2$), all eigenvalues are non-negative, hence $\gamma_*=1-\lambda_2$, whereas for periodic chains $t_{\rm rel}=+\infty$. 

The following result relates the mixing time of a reversible and irreducible Markov chain to the relaxation time. 

\begin{lemma}[Theorem~12.4 of \cite{levin2017markov}]\label{lem:mix_rel}
Let $P$ the transition matrix of a reversible irreducible Markov chain with state space $S$ and $\pi_{\rm min}=\min_{x\in S}\pi(x)$. Then 
\[
t_{\rm mix}(\varepsilon)\leq t_{\rm rel}\log\left(\frac{1}{\varepsilon\pi_{\rm min}}\right). 
\]
\end{lemma}

The second result that we will use relates the average hitting times of states of a reversible irreducible Markov chain to the eigenvalues of the corresponding transition matrix. (Note that in \cite{aldous-fill-2014} this result is referred to as ``the eigentime identity.'')

\begin{lemma}[Theorem~15 of \cite{Broder1989}; see also Proposition~3.13 of \cite{aldous-fill-2014} or Lemma 12.17 \cite{levin2017markov}]\label{lem:eigentime}
For a reversible and irreducible Markov chain on a finite state space $S$ it holds that 
\[
\sum_{x\in S}\sum_{y \in S}\pi(x)\pi(y)E_{x}[H_y] =\sum_{y \in S}\pi(y)E_{x}[H_y]= \sum_{i=2}^{|S|}\frac{1}{1-\lambda_i}, 
\]
where $H_y$ is the hitting time of $y$ and the first equality is a consequence of the ``random target Lemma'' (see, e.g., \cite[Lemma 2.29]{aldous-fill-2014} or \cite[Lemma 10.1]{levin2017markov}).
\end{lemma}

The last auxiliary result, which is taken from \cite{diaconis1991},  provides a lower bound on the smallest eigenvalue $\lambda_{|S|}$ of a reversible, irreducible and aperiodic Markov chain in terms of geometric quantities of the corresponding graph. Before stating the result we need to introduce some auxiliary concepts.
Let us consider the graph with vertex set $S$ and an edge from $x$ to $y$ if $P(x,y)>0$. 
For $x \in S$, let $\sigma_x$ be a path from $x$ to $x$ \emph{with an odd  number of edges} (such a path always exists for irreducible and aperiodic chains) and let $\Sigma$ be a collection of such paths, one for each $x \in S$. Let $|\sigma_x|_P$ be the weighted (by resistance) length of the  path $\sigma_x$, defined as 
\[
|\sigma_x|_P:=\sum_{\{z,y\} \in \sigma_x}\frac{1}{\pi(z)P(z,y)}.
\]
A geometric quantity of interest to us is given by
\[
i(\Sigma):= \max_{\{z,y\}}\sum_{\sigma_x \ni \{z,y\}}|\sigma_x|_P \, \pi(x).
\]

\begin{lemma}[Proposition~2 of \cite{diaconis1991}]\label{lem:geom_eigen}
For a reversible, irreducible and aperiodic Markov chain with finite state space $S$ and transition matrix $P$, the smallest eigenvalue $\lambda_{|S|}$ satisfies
\[
\lambda_{|S|}\geq -1 + \frac{2}{i(\Sigma)}.
\]
\end{lemma}

\begin{remark}\label{rem:geom_eigen}
Note that for random walk on a rooted tree $T=(S,E)$ with a unique self-loop at the root, for each vertex $x$ we can take, as $\sigma_x$, the shortest path from $x$ to the root, through the loop, then back to $x$, obtaining $|\sigma_x|_P=2|E|(2\, \mathsf{ dist}_T(x,{\rm root})+1)$ and
\[
i(\Sigma)= \sum_{x\in S} \left(2\, \mathsf{dist}_T(x,{\rm root})+1\right)\mathsf{deg}_{T}(x)\leq \left(2\, \mathsf{diam}(T) + 1\right)2|E|.
\]
Thus, from Lemma~\ref{lem:geom_eigen},  the corresponding smallest eigenvalue  satisfies $\lambda_{|S|}\geq -1 + \left((2\mathsf{diam}(T)+1)|E|\right)^{-1}$.
\end{remark}

We now have  all the ingredients to provide an upper bound on the relaxation time.

\begin{theorem}[Relaxation time bound]
For simple random walk on a tree $T=(S,E)$ with a (unique) self-loop, it holds that
\[
t_{\rm rel}\leq (2\,\mathsf{ diam}(T)+1)|E|.
\]
\end{theorem}

\begin{proof}
There are two possible scenarios: $\max\{|\lambda_2|,|\lambda_{|S|}|\}=|\lambda_2|$ or $\max\{|\lambda_2|,|\lambda_{|S|}|\}=|\lambda_{|S|}|$. In the first case, necessarily $\lambda_2\ge 0$, and thus by Lemma~\ref{lem:eigentime} we immediately get that $t_{\rm rel}\leq \sum_{y \in S}\pi(y)E_{x}[H_y]\leq \max_{x,y}E_{x}[H_y]\leq |E| \, \mathsf{diam}(T)$, where the last inequality follows from the commute time identity \cite[Proposition 10.7]{levin2017markov}. In the second case, $\lambda_{|S|}\leq 0$ and by Lemma~\ref{lem:geom_eigen} and Remark~\ref{rem:geom_eigen}, we obtain that  $\gamma*=1-|\lambda_{|S|}|\geq \left((2\,\mathsf{diam}(T)+1)|E|\right)^{-1}$, and thus the claim follows.
\end{proof}

Using Lemma~\ref{lem:mix_rel}, we immediately obtain the following corollary. 

\begin{corollary}[Mixing time bound]\label{cor:mixingbound}
For simple random walk on a tree $T=(S,E)$ with a (unique) self-loop, it holds that
\[
t_{\rm mix}(\varepsilon) \leq  (2\,\mathsf{ diam}(T)+1)|E| \log \left(\frac{2|E|}{\varepsilon}\right).
\]
\end{corollary}

\begin{remark}
Y.~Peres outlined a proof (via personal correspondence) that gives a slightly stronger bound: namely, the logarithmic factor can be dropped. The argument uses results on averaging distributions over two consecutive steps \cite{HerPesPower}. For our purposes, the presence of the logarithmic factor is rather irrelevant though.
\end{remark}

\subsection{Strong Stationary Times}\label{ss:stt}

The coupling between the BA-model and the TBRW relies on the framework of {\it strong stationary times} for Markov chains introduced in \cite{aldous1986shuffling,aldous1987strong}. In this subsection we will introduce the main concepts and  results we will need regarding stationary times and mixing times. The reader who is already familiar with these notions, may skip this and proceed  directly to Lemma \ref{lem:remainstationary}.

\begin{definition}[Separation distance]\label{def:sx} Given a Markov chain on a finite state space $S$,  transition matrix $P$  and stationary distribution $\pi$, and $x\in S$, we call the separation distance from $x$ the following quantity:
\begin{equation}\label{eq:sx}
    s_x(t) := \max_{y\in S}\left[1-\frac{P^t(x,y)}{\pi(y)}\right].
\end{equation}
We also introduce the separation distance
\begin{equation}\label{eq:s}
    s(t) := \max_{x\in S}s_x(t).
\end{equation}
\end{definition}

\begin{definition}[Strong stationary time] Let $\{Y_t\}_{t\ge 0}$ be a Markov chain with stationary distribution $\pi$, and adapted to a filtration $\{\mathcal{F}_t\}_{t\ge 0}$. A strong stationary time for $\{Y_t\}_t$ and starting point $y_0$ is a $\{\mathcal{F}_t\}_t$-stopping time $\eta$ such that $Y_\eta$ is distributed as $\pi$ and it is independent of $\eta$. In symbols, for any state $y$ and $k \in \mathbb{N}$, we have
\begin{equation}
    P_{y_0} \left( X_\eta = y,  \eta=k \right) = \pi(y)P_{y_0}\left(\eta = k\right).
\end{equation}
\end{definition}
The connection between strong stationary times and separation distance is given by the following result.

\begin{proposition}[Proposition 3.2 of \cite{aldous1987strong}]\label{prop:optimalstrong} Let $\{Y_t\}_t$ be a Markov chain with initial state $y_0$. Then, there exists a strong stationary time $\eta_{y_0}$, such that for all $t$,
\begin{equation}
    P_{y_0}(\eta_{y_0} > t)=s_{y_0}(t).
\end{equation}
We will call $\eta_{y_0}$ an {\bf optimal} strong stationary time. (The reason for the name is that $P_{y_0}(\eta_{y_0} > t) \ge s_{y_{0}}(t)$ holds for \emph{any} strong stationary time.) \end{proposition}

It will be useful for us to relate separation distance and total variation distance. Under the same settings as Definition~\ref{def:sx}, we let
\begin{equation}\label{eq:dt}
    d(t):= \max_{x\in S} \| P^t(x, \cdot) - \pi \|_{TV}.
\end{equation}
As usual, we set $t_{mix}:=\min\{t : d(t) \leq 1/4\}$.

\begin{proposition}[Lemma 6.17 of \cite{levin2017markov}]\label{prop:sepdtv} For any reversible Markov chain,
$$
s(2t) \le 4d(t).
$$
\end{proposition}

The following result relates the tail of an optimal strong stationary time to $t_{mix}$.

\begin{proposition}\label{prop:ssttail}Let $\{Y_t\}_t$ be an aperiodic and irreducible Markov chain. Also, let $\eta_{x_0}$ be an optimal strong stationary time for the chain started at $x_0$. Then, for any $\ell \in \mathbb{N}$,
\begin{equation*}
    P_{x_0}\left( \eta_{x_0} > \ell t_{mix} \right) \le 4\cdot 2^{-\ell/2}.
\end{equation*}
\end{proposition}
\begin{proof} The proof follows by combining Proposition \ref{prop:optimalstrong} and Proposition \ref{prop:sepdtv} with the fact that for any natural number $\ell$ 
$$
d(\ell t_{mix}) \le 2^{-\ell};
$$
see \cite[Section 4.5]{levin2017markov}. Then,
$$
P_{x_0}(\eta_{x_0} > \ell t_{mix}) \le 4d(\ell t_{mix}/2) \le 4 \cdot 2^{-\ell/2},
$$
as desired.
\end{proof}

If $\eta$ is a stationary time, then all the subsequent deterministic steps of the random walk will be distributed according to the stationary distribution. However, this is not necessarily the case for random times that come after the stationary time. As a matter of fact, if one considers $\tau_v$,  the hitting time of a vertex $v$, then  $X_{\tau_v} \equiv v$. Nevertheless, the next lemma tells us that whenever a random time $\tau$ is independent of $X$ and $\eta$, then conditioned on $\{\eta<\tau\}$ we have $X_{\tau-1} \sim \pi$. 

\begin{lemma}[Remaining stationary after a random amount of steps]\label{lem:remainstationary} Let $\{X_n\}_{n\in \mathbb{N}}$ be a SSRW (Simple Symmetric Random Walk) on a rooted tree $T$ with a self-loop at the root. Also, let $\pi$ be its stationary distribution, and $\eta$ a strong stationary time given by Proposition~\ref{prop:optimalstrong}. Additionally, let $\tau$ be a random time independent of $X$ and $\eta$, satisfying that $0<\tau<\infty$ a.s. Then,
$$
    P_{v_0} \left(X_{\tau-1} = v \; \middle | \; \eta < \tau \right) = \pi(v).
$$
\end{lemma}

\begin{proof} The claim follows from the (strong) Markov property. In detail:
\begin{equation*}
    \begin{split}
    P_{v_0} \left(X_{\tau-1} = v \; \middle | \; \eta < \tau \right) & = \sum_{k < \ell} P_{v_0} \left(X_{\ell-1} = v \; \middle | \; \eta=k,\ \tau=\ell \right) \frac{ P_{v_0} ( \eta=k,\ \tau=\ell ) }{P_{v_0} (\eta < \tau )}\\
     & = \sum_{k < \ell} P_{v_0} \left(X_{\ell-1} = v \; \middle | \; \eta=k \right) \frac{  P_{v_0} ( \eta=k,\ \tau=\ell ) }{P_{v_0} (\eta < \tau )}\\
     & = \sum_{k < \ell} \sum_w P_{v_0} \left(X_{\ell-1} = v \; \middle | \; X_k=w \right) P_{v_0} \left( X_k=w\; \middle | \; \eta=k \right)\\
     & \quad\times\, \frac{  P_{v_0} ( \eta=k,\ \tau=\ell ) } {P_{v_0} (\eta < \tau )}\\
      & = \sum_{k < \ell} \sum_w P_{v_0} \left(X_{\ell-1} = v \; \middle | \; X_k=w \right) \pi(w) \frac{  P_{v_0} ( \eta=k,\ \tau=\ell ) }{P_{v_0} (\eta < \tau )}\\
      & = \sum_{k < \ell} \pi(v) \frac{ P_{v_0} ( \eta=k,\ \tau=\ell ) }{P_{v_0} (\eta < \tau )} = \pi(v),
    \end{split}
\end{equation*}
where, to get the second line we used that $\tau$ is independent of everything, to get the third line we used the Markov property, to get the fourth line we used that $\eta$ is a strong stationary time, and to get the fifth line we used that $\pi$ is stationary.  
\end{proof}

\subsection{Strong Stationary times for the TBRW}\label{ss:ssttbrw}

The coupling given in the next section will rely on the existence of optimal strong stationary times discussed in the previous subsection. However, it is important to better explain what we mean by strong stationary times in the context of the TBRW, because the tree over which the walker is walking changes over time. Observe that, given an initial state $(T,x)$, the walker of the TBRW with a sequence of laws $\mathcal{L}$, until the appearance of new vertices, can be coupled to agree with a SSRW on $T$ starting at $x$. After that, we can just let the two walkers evolve independently of each other. With this coupling in mind, it makes sense to talk about stationary times for the TBRW.

More precisely, consider a TBRW $X$ starting at $(T,x)$, and let $\widetilde{X}$ be a SSRW on $T$ started at $x$. Also, let $\tau$ be the stopping time when the walker adds new vertice(s) for the first time; that is,
\begin{equation}\label{def:tau}
    \tau = \inf \{n>0 \, : \, Z_n \ge 1\}.
\end{equation}
Now let $P_{T,x; \mathcal{L}}$ be the coupling measure of $X$ and $\widetilde{X}$ that we mentioned above: $X_i = \widetilde{X}_i$ for $i\leq \tau$, while $\{X\}_{i\ge \tau}$ and $\{\widetilde{X}\}_{i\ge \tau}$, given $\tau$ and the vertex $X_\tau=\widetilde{X}_\tau$, are independent. Furthermore, let $\eta_x$ be a strong stationary time for $\widetilde{X}$, and let $\pi_T$ be the stationary distribution of $\widetilde{X}$:
\begin{equation}\label{def:piT}
    \pi_T(v) := \frac{\degree{T}{v}}{\sum_{z \in T}\degree{T}{z}}.
\end{equation} 
The next result concerns the distribution of $X_{\eta_x}$ in the coupling $P_{T,x; \mathcal{L}}$.

\begin{proposition}\label{prop:rightbias}Let $T$ be a finite tree, $x, v$ two of its vertices, and consider a sequence of laws $\mathcal{L}$ satisfying that $\mathcal{L}_n(\{0\})>0$ for all $n$. Then, conditioned on $\{\eta_x < \tau\}$, both $X_{\eta_x}$ and $X_{\tau -1}$ follow the distribution $\pi_T$. In symbols,
\begin{equation}
    P_{T,x;\mathcal{L}}\left( X_{\eta_x} = v \; \middle | \; \eta_x < \tau \right) = P_{T,x;\mathcal{L}}\left( X_{\tau-1} = v \; \middle | \; \eta_x < \tau \right) = \pi_T(v).
\end{equation}
\end{proposition}

We highlight that the need for the sequence of laws to satisfy $\mathcal{L}_n(\{0\})>0$ for all $n$ is to make sure that the event $\{\eta_x < \tau\}$ has positive probability.

\begin{proof} Using that $\eta_x$ is a strong stationary time for $\widetilde{X}$ and is independent of $\tau$, we have that 
\begin{equation}
    \begin{split}
        P_{T,x; \mathcal{L}}\left( X_{\eta_x} = v, \eta_x < \tau \right) & = P_{T,x; \mathcal{L}}\left( \widetilde{X}_{\eta_x} = v, \eta_x < \tau \right)\\
        &= \sum_{n=0}^\infty P_{T,x; \mathcal{L}}\left( \widetilde{X}_{\eta_x} = v, \eta_x = n\right)P_{T,x; \mathcal{L}}\left( n < \tau \right)\\
        &= \pi_T(v)\sum_{n=0}^\infty P_{T,x; \mathcal{L}}\left(\eta_x = n\right)P_{T,x; \mathcal{L}}\left( n < \tau \right)\\
        & = \pi_T(v) P_{T,x; \mathcal{L}}\left(\eta_x < \tau \right).
    \end{split}
\end{equation}
Also notice that from the perspective of $\widetilde{X}$, $\tau$ is just a random clock according to an independent source of randomness. Thus, using that both walkers are the same up to time $\tau -1$ and Lemma \ref{lem:remainstationary}, we also have that for any vertex $v \in T$,
\begin{equation}
\begin{split}
P_{T,x; \mathcal{L}} \left( X_{\tau-1}=v, \eta_x < \tau \right) & = P_{T,x; \mathcal{L}} \left( \widetilde{X}_{\tau-1}=v, \eta_x < \tau \right) = \pi_T(v)P_{T,x; \mathcal{L}} \left( \eta_x < \tau \right),
\end{split}
\end{equation}
which implies that conditioned on $\{\eta_x < \tau \}$, $X_{\tau -1}$ is distributed according to $\pi_T$.
\end{proof}

This proposition will allow us to talk about stationary times for the TBRW. Thus, from now on, whenever we say $\eta_x$ is a strong stationary time for the TBRW starting at $(T,x)$ and with sequence of laws $\mathcal{L}$, the reader must have in mind the construction discussed above.

\section{Proof of Theorem \ref{thm:generaltransfer}: General Condition for Transfer Principle}\label{sec:generaltransfer}
In this section, we will work with the general condition \mixc\ (see Definition \ref{def:TP}), which, once being satisfied by the \name, gives us a way to transfer results from the BA-model to the \name\ (Theorem~\ref{thm:generaltransfer}). The proof relies on a coupling argument with the BA-model done in Section~\ref{ss:coupling}. Thus, we will first construct such a coupling and then, in Section~\ref{ss:proofgeneraltransfer}, we show how to obtain the transfer principle from it.

\subsection{Coupling with the preferential attachment (BA) model}\label{ss:coupling} 

In this section we will construct a coupling between the tree sequence produced by the TBRW and the one produced by the BA model. 

 \begin{theorem}[Coupling with preferential attachment]\label{thm:coupling.w.PA} 
 Consider a TBRW satisfying condition \mixc. Then, there exists a probability space containing the TBRW and a sequence of random graph processes $\{\boldsymbol{\mathsf{G}}_n\}_{n \in \mathbb{N}}$, where $\boldsymbol{\mathsf{G}}_n:= \{G_k^{(n)}\}_{k \in \mathbb{N}}$ is a BA-tree starting at $T_{\tau_n}$ such that 
 $$
\lim_{n \to \infty} D_{TV}\big( \{T_{\tau_{n+k}}\}_{k\in \mathbb{N}}, \{G^{(n)}_k\}_{k\in \mathbb{N}} \big) = 0\,. 
 $$
\end{theorem}

\begin{proof} In the coupling, we will have a single TBRW process $(T_k,X_k)_{k\ge 0}$, with growth times $(\tau_n)_{n\ge 0}$. We will also use a sequence $(\eta_n)_{n\ge 0}$ of optimal strong stationary times for SSRW random walk on $T_{\tau_n}$ started at $X_{\tau_n}$, independently for different $n$'s, which exist by Proposition~\ref{prop:optimalstrong} and Subsection~\ref{ss:ssttbrw}. 

For each $n\ge 0$, we will now construct $\{G_k^{(n)}\}_{k \in \mathbb{N}}$ distributed as a BA process starting from $G^{(n)}_0 =T_{\tau_n}$, inductively in $k$. 

Enlarging our probability space if necessary, assume that for each $n$ we have a sequence of independent random variables $\{U^{(n)}_k\}_k$, which is also independent of the TBRW process and $U^{(n)}_k \sim \mathsf{Uni}[0,1]$ for all $n$ and $k$. Then, we generate $G^{(n)}_1$ as follows:
\begin{enumerate}
    \item[(A)] If $\eta_n  <  \tau_{n+1} - \tau_n =: \Delta\tau_{n}$, then we set $G^{(n)}_1 = T_{\tau_{n+1}}$;
    \item[(B)] Otherwise, we use (in a way detailed below) the independent source of randomness provided by $U^{(n)}_1$, and independently of everything else, to select a vertex of $v\in T_{\tau_n}$ with probability 
    $$
    \frac{\mathsf{deg}_{T_{\tau_n}}(v)}{\sum_{z\in T_{\tau_n}}\mathsf{deg}_{T_{\tau_n}}(z)},
    $$
    and connect it to a new vertex add to $G_0^{(n)} := T_{\tau_n}$. This new graph is then $G^{(n)}_1$. 
\end{enumerate}
Let us already point out here that option (A) will be the ``useful one'' for us: if that occurs, then we will be able to continue the coupling by generating $G^{(n)}_2$ starting from $G^{(n)}_1 = T_{\tau_{n+1}}$. On the other hand, if (B) occurs, then we cannot relate $G_1^{(n)}$ to $T_{\tau_{n+1}}$, hence we consider the coupling to be a failure, and will just generate the futures independently. 

Regarding the use of the uniform random variable at (B). This can be done as follows. We first enumerate the vertices of $T_{\tau_n}$. Then, we split the interval $[0,1]$ into $|T_{\tau_n}|$ subintervals in a way that the length of the $j$-th subinterval is exactly
$$
    \frac{\mathsf{deg}_{T_{\tau_n}}(v_j)}{\sum_{z\in T_{\tau_n}}\mathsf{deg}_{T_{\tau_n}}(z)}.
$$
Then, if $U^{(n)}_1$ belongs to the $j$-th subinterval, we connect the new vertex to $v_j$.

Our first claim is that $G^{(n)}_1$ is distributed as the first step of a BA random graph started from $G^{(n)}_0=T_{\tau_n}$. Indeed, let $\{\mathcal{F}_n\}_n$ be the canonical filtration for the \name\ and let's call a graph $H$ admissible for $T_{\tau_n}$ if $H$ is $T_{\tau_n}$ together with a new vertex $v_*$ connected to some vertex $u$ of $T_{\tau_n}$. Then, observe that
\begin{equation}\label{eq:g1}
    \begin{split}
        P_{T_0,x_0; \laws^{}}\left( G^{(n)}_1 = H  \middle | \mathcal{F}_{\tau_n}\right) & = P_{T_0,x_0; \laws^{}}\left( G^{(n)}_1 = H, \eta_n  <  \Delta\tau_n \middle | \mathcal{F}_{\tau_n} \right) \\
        &\quad +  P_{T_0,x_0; \laws^{}}\left( G^{(n)}_1 = H, \eta_n \geq \Delta \tau_n  \middle | \mathcal{F}_{\tau_n} \right).
    \end{split}
\end{equation}
Notice that when $\eta_n  \geq \Delta\tau_n$, we generate $G^{(n)}_1$ using the independent source of randomness coming from $U^{(n)}_1$, thus using independence of $U^{(n)}_1$
\begin{equation}\label{eq:g12}
    \begin{split}
        P_{T_0,x_0; \laws^{}}\left( G^{(n)}_1 = H,  \eta_n \geq  \Delta \tau_n \middle | \mathcal{F}_{\tau_n}\right) = \frac{\mathsf{deg}_{T_{\tau_n}}(u)}{\sum_{z\in T_{\tau_n}}\mathsf{deg}_{T_{\tau_n}}(z)}P_{T_0,x_0; \laws^{}}\left(  \eta_n \geq \Delta \tau_n \middle | \mathcal{F}_{\tau_n}\right).
    \end{split}
\end{equation}
For the first term of the RHS of \eqref{eq:g1}, observe that if $v_*$ denotes the vertex added at time $\tau_{n+1}$, then $v_*$ connects to $u$ if, and only if, $X_{\tau_{n+1} -1}=u$. Thus, by the strong Markov property, 
\begin{equation}\label{eq:g13}
    \begin{split}
         P_{T_0,x_0; \laws^{}}&\left( G^{(n)}_1 = H, \eta_n < \Delta \tau_n \ \middle |\ \mathcal{F}_{\tau_n}\right)  
         =  P_{T_0,x_0; \laws^{}}\left( X_{\tau_{n+1} -1} = u, \eta_n < \Delta \tau_n \mid\mathcal{F}_{\tau_n}\right)
         \\
         &= P_{T_{\tau_n},X_{\tau_n}; \laws^{(\tau_n)}} \left(  X_{\tau_{1} -1} = u, \eta_0 <  \tau_1 \right) \\
         &= \frac{\mathsf{deg}_{T_{\tau_n}}(u)}{\sum_{z\in T_{\tau_n}}\mathsf{deg}_{T_{\tau_n}}(z)} P_{T_{\tau_n},X_{\tau_n}; \laws^{(\tau_n)}} \left( \eta_0 <  \tau_1 \right)\\
         &= \frac{\mathsf{deg}_{T_{\tau_n}}(u)}{\sum_{z\in T_{\tau_n}}\mathsf{deg}_{T_{\tau_n}}(z)} P_{T_0,x_0; \laws^{}}\left(  \eta_n < \Delta \tau_n \mid \mathcal{F}_{\tau_n}\right) ,
    \end{split}
\end{equation}
where in the third equality we used Proposition~\ref{prop:rightbias}. Plugging \eqref{eq:g12} and \eqref{eq:g13} into~\eqref{eq:g1} gives us 
$$
    P_{T_0,x_0; \laws^{}}\left( G^{(n)}_1 = H\ \middle |\ \mathcal{F}_{\tau_n}\right) = \frac{\mathsf{deg}_{T_{\tau_n}}(u)}{\sum_{z\in T_{\tau_n}}\mathsf{deg}_{T_{\tau_n}}(z)},
$$
which means that under $P_{T_0,x_0; \laws^{}}$, the random graph $G^{(n)}_1$ is distributed as one step of a BA random graph started from the random initial condition $T_{\tau_n}$.

Now the strong Markov property of $\{(T_n,X_n)\}_n$ allows us to iterate the above procedure to generate $G^{(n)}_2$ and the further steps of the BA model. Indeed, assume we have successfully constructed $\{G^{(n)}_i\}_{i=0}^k$, which represents $k$ steps of the BA model with $G_0^{(n)} = T_{\tau_n}$. Then, in order to generate~$G^{(n)}_{k+1}$ we do as before:
\begin{enumerate}
    \item[(A)] If 
    $\eta_{n+i} < \Delta\tau_{n+i}$, for all $i\le k$, then we set $G^{(n)}_{k+1} = T_{\tau_{n+k+1}}$;
    \item[(B)] Otherwise, we use the independent source of randomness provided by $U^{(n)}_{k+1}$, and independently of everything else select a vertex of $v\in G^{(n)}_k$ with probability 
    $$
    \frac{\mathsf{deg}_{G^{(n)}_k}(v)}{\sum_{z\in G^{(n)}_k}\mathsf{deg}_{G^{(n)}_k}(z)},
    $$
    and connect it to a new vertex added to $G^{(n)}_k$. This new graph is then $G^{(n)}_{k+1}$. 
\end{enumerate}

We claim that $G^{(n)}_{k+1}$ is distributed as the $k+1$-th step of the BA-model started from the random graph $T_{\tau_n}$. The proof is by induction on $k$, and is exactly as before, just the partition~\eqref{eq:g1} is now according to whether the event
\begin{equation}\label{eq:Ak}
    A^{(n)}_{k} = \bigcap_{i=0}^{k}\left \lbrace \eta_{n+i}<\Delta \tau_{n+i}\right 
    \rbrace
\end{equation}
occurs or not. If yes, we follow the argument of ~\eqref{eq:g13}; if not, then follow the argument of~\eqref{eq:g12}, using that $U_{k+1}^{(n)}$ is independent of everything else.

Now notice that for any $n$, we have the following inclusion of events
\begin{equation}
\bigcap_{i=1}^\infty\left \lbrace \eta_{n+i} < \Delta \tau_{n+i}\right \rbrace \subset \left \lbrace T_{\tau_{n+i}} = G^{(n)}_i, \forall i \right \rbrace,
\end{equation}
which implies that
\begin{equation*}
    D_{TV}(\{T_{\tau_{n+k}}\}_{k\in \mathbb{N}}, \{G^{(n)}_k\}_{k\in \mathbb{N}}) \le P_{T_0,x_0; \mathcal{L}^{}}\left( \bigcup_{i=1}^\infty\left \lbrace \eta_{n+i} \geq  \Delta \tau_{n+i}\right \rbrace \right). 
\end{equation*}
Finally, by the definition of $\limsup$ of events and part (2) of condition \mixc{} we have that
$$
\lim_{n \to \infty} D_{TV}\big(\{T_{\tau_{k+n}}\}_{k\in \mathbb{N}}, \{G^{(n)}_k\}_{k\in \mathbb{N}}\big) \le P_{T,x;\mathcal{L}}\left(\eta_k \geq \Delta\tau_k \text{, i.o.} \right) = 0,
$$
which concludes the proof.
\end{proof}

\subsection{Proof of Theorem \ref{thm:generaltransfer}: The Transfer Principle under \mixc{}}\label{ss:proofgeneraltransfer}
We can now prove that under condition \mixc{} we can transfer properties of the BA-tree to the \name.
\begin{proof}[Proof of Theorem \ref{thm:generaltransfer}] Fix an initial condition for the \name\ $(T,x)$. Also let $M$ be the following random variable
\begin{equation}\label{def:M}
    M := \inf \{ k \, : \, \eta_n < \Delta \tau_n \text{ for all }n \ge k \},
\end{equation}
where $\eta_{n}$ is the optimal stationary time for the TBRW started at~$(T_{\tau_{n}},X_{\tau_{n}})$ and with sequence of laws $\mathcal{L}^{(\tau_{n})}$, independently for different $n$'s. Notice that, by item (2) of condition \mixc, $M$ is finite almost surely. Moreover, recall the proof of Theorem \ref{thm:coupling.w.PA}, where we constructed a probability space containing a sequence of BA-trees. For each $n$, $\{G^{(n)}_k\}_{k \in \mathbb{N}}$ is a BA-model starting from~$T_{\tau_n}$.  Also recall from the construction that, on the event $\{M= j\}$, 
$$
T_{\tau_{j+k}} = G^{(j)}_k, \text{ for all } k\in \mathbb{N}.
$$
Now, let $\mathcal{G}$ be an asymptotic graph property for the BA-tree that holds $P_{G_0^{\rm loop}}$-almost surely.  Then,  denoting by $Q_{T_{\tau_j}}$  the distribution of $\{G^{(j)}_k\}_{k \in \mathbb{N}}$,  necessarily,
\begin{equation}
    Q_{T_{\tau_j}}\left(\{G^{(j)}_k\}_{k \in \mathbb{N}} \in \mathcal{G} \right) \equiv 1, \qquad 
    P_{T,x; \mathcal{L}^{}}\text{-a.s.}
\end{equation}
Thus, for all $j \in \mathbb{N}$,
\begin{equation}
    \begin{split}
        P_{T,x; \mathcal{L}^{}}\left( \{T_{\tau_{j+k}}\}_k \in \mathcal{G}, M=j\right) & = P_{T,x; \mathcal{L}^{}}\left( \{G^{(j)}_{k}\}_k \in \mathcal{G}, M=j\right) \\
        & = P_{T,x; \mathcal{L}^{}}\left( M=j\right);
    \end{split}
\end{equation}
in the second equality, we used the fact that the event $\{\{G^{(j)}_{k}\}_k \in \mathcal{G}\}$ has total probability. By the above and using the hypothesis that $\mathcal{G}$ is an asymptotic graph property, we have that 
$$
P_{T,x; \mathcal{L}^{}}\left( \{T_{\tau_{k}}\}_k \in \mathcal{G}, M=j\right) = P_{T,x; \mathcal{L}^{}}\left( \{T_{\tau_{j+k}}\}_k \in \mathcal{G}, M=j\right) = P_{T,x; \mathcal{L}^{}}\left( M=j\right).
$$
Finally, summing over $j$ and recalling that $M$ is finite almost surely gives the desired result.
\end{proof}

\section{Proof of Theorem \ref{thm:transfer}: An example satisfying condition \mixc}\label{sec:transfer}

In the rest of the paper, we will always consider a \name\ with a sequence of laws $\mathcal{L}_n = {\rm Ber}(p_n)$, where $p_n$ satisfies the following condition:

\begin{definition}\label{def:G}
We say that a sequence of probabilities $(p_n)_{n\ge 1}$ satisfies condition (G) if it is decreasing, and there exist absolute constants $0<\g,\G<\infty$ and $\gamma \in [0,1]$ such that 
$\g n^{-\gamma} < p_n < \G n^{-\gamma}$ holds for all $n\ge 1$.
\end{definition}

In this section, we will prove Theorem~\ref{thm:transfer}, which will follow from the two auxiliary results below.

\begin{theorem}[Upper bound on the diameter]\label{thm:diam} 
Consider an $\mathcal{L}$-\name\ where $\mathcal{L}_n = \mathsf{Ber}(p_n)$ with $p_n$ satisfying Definition~\ref{def:G} with $\gamma \in (1/2,1]$.  Then, for any $\varepsilon >0$,  $M>0$, and and initial finite tree $T$, there exist positive constants $C$ and $C'$ depending on $\varepsilon$, $M$, $\gamma$, $\g$, $\G$ and $|T|$, such that for all $n\geq 1$
$$
\mathbb{P}_{T,x;\mathcal{L}^{}} \left( \mathsf{diam}(T_n) \ge Cn^{\varepsilon}\right) \le \frac{C'}{n^M}.
$$
\end{theorem}

As it will become clear later, Theorem~\ref{thm:diam} is a key step in our bootstrap argument. It provides a ``good enough'' bound that enables us to show that the walker will mix before adding new leaves. However, it is also important to point out that, when condition (M) is satisfied,  the order given in Theorem~\ref{thm:diam} is far from the correct order of the diameter of $T_n$, which is of order $\log n$ (see Corollary~\ref{cor:height}).

For the sake of readability, the proof of the above result will be postponed to Section~\ref{sec:diameter}.

\begin{lemma}[Growth events do not occur too early] Consider a TBRW starting at $(T,x)$ with sequence of laws $\mathcal{L}_n \sim \mathsf{Ber}(p_n)$ with $p_n$ satisfying Definition~\ref{def:G} with  $\gamma \in (2/3,1)$. Let $\tau_k$ be the growth times as in \eqref{def:tauk}. Then, for any $\varepsilon < \frac{\gamma}{1-\gamma} - 2$
$$
\Pd_{T,x; \mathcal{L}^{}}\left( \Delta \tau_i < i^{1+\varepsilon}, \text{ i.o.} \right) = 0.
$$
\end{lemma}

\begin{proof} Recall that  $Z_k$ is the indicator that a leaf has been added at time $k$. Also fix $\delta>0$  to be chosen properly later. By the definition of $\tau_n$ and a Chernoff bound \cite[Theorem 2.7]{mcdiarmid1998concentration}, we have that there exists a positive constant $C$ depending on  $\gamma$ only, such that
\begin{equation}\label{eq:taunprob}
    \Pd_{T,x; \mathcal{L}^{}}\left( \tau_n < n^{\frac{1-\delta}{1-\gamma}} \right) \leq \Pd_{T,x; \mathcal{L}^{}}\left( \sum_{k=1}^{n^{\frac{1-\delta}{1-\gamma}}} Z_k \ge n \right) \le e^{-Cn},
\end{equation}
using that $$
\mathsf{E}\left(\sum_{k=1}^{n^{\frac{1-\delta}{1-\gamma}}} Z_k\right) \asymp \mathsf{Var}\left(\sum_{k=1}^{n^{\frac{1-\delta}{1-\gamma}}} Z_k\right) \asymp n^{1-\delta}.
$$
On the other hand, by the strong Markov property, it follows that
\begin{equation}
    \Pd_{T,x; \mathcal{L}^{}}\left( \Delta \tau_n \ge n^{1+\varepsilon} \mid \mathcal{F}_{\tau_{n-1}}\right) \ge \prod_{k=1}^{n^{1+\varepsilon}}\left( 1- \frac{\G}{(\tau_{n-1}+k)^{\gamma}}\right),
\end{equation}
which implies that, on the event $\{\tau_{n-1} \ge (n-1)^{(1-\delta)/(1-\gamma)}\}$, we have
\begin{equation}
    \Pd_{T,x; \mathcal{L}^{}}\left( \Delta \tau_n \ge n^{1+\varepsilon} \mid \mathcal{F}_{\tau_{n-1}}\right) \ge \left( 1- \frac{\G}{(n-1)^{(1-\delta)\gamma/(1-\gamma)}}\right)^{n^{1+\varepsilon}}.
\end{equation}
We now choose $\delta>0$ small enough so that $(1-\delta)\gamma/(1-\gamma)-(1+\varepsilon) > 1$. Observing that, by \eqref{eq:taunprob}, the event $\{\tau_{n-1} \ge (n-1)^{(1-\delta)/(1-\gamma)}\}$ occurs with probability at least $1-e^{-Cn}$, we deduce that there exists a $\delta'>0$ (depending on $\varepsilon$ and $\gamma$) such that
\begin{equation}\label{ineq:deltatauk}
    \Pd_{T,x; \mathcal{L}^{}}\left( \Delta \tau_n \le n^{1+\varepsilon} \right) \le \frac{1}{n^{1+\delta'}} . 
\end{equation}
This, combined with the Borel-Cantelli lemma, gives  us the result.
\end{proof}

Now we have all the results we need, we proceed to the proof of Theorem \ref{thm:transfer}.

\begin{proof}[Proof of Theorem \ref{thm:transfer}]
In light of Theorem \ref{thm:generaltransfer}, it is enough show that the \name\ with the sequence of laws $\mathcal{L}_n = {\rm Ber}(p_n)$ satisfying condition (G) of Definition~\ref{def:G} with $\gamma \in (2/3,1]$ also satisfies condition \mixc. Notice that item (1) of \mixc\ is clearly satisfied. Then, we need to show item (2) of \mixc, meaning that
\begin{equation}\label{eq:tp2}
P_{T,x;\mathcal{L}^{}}\left(\eta_k \geq \Delta\tau_k \text{, i.o.} \right) = 0,
\end{equation}
holds for any finite initial condition $(T,x)$ and time shift $m \in \mathbb{N}$.
 We will deal with the case $\gamma \in (2/3,1)$ first, which is the harder one.\\

\noindent {\it \underline{Case $\gamma \in (2/3,1)$.} }
Observe that for a fixed $k$ and sufficiently small $\varepsilon>0$, \eqref{ineq:deltatauk} yields to
\begin{equation}\label{ineq:first}
    \begin{split}
        P_{T,x;\mathcal{L}^{}}\left(\eta_k \geq \Delta\tau_k  \right) & \le P_{T,x;\mathcal{L}^{}}\left(\eta_k > k^{1+\varepsilon} \right) + \Pd_{T,x;\mathcal{L}^{}}\left(\Delta \tau_k \le k^{1+\varepsilon}\right)\\
        & \le P_{T,x;\mathcal{L}^{}}\left(\eta_k > k^{1+\varepsilon} \right) + \frac{1}{k^{1+\delta'}},
    \end{split}
\end{equation}
for some positive $\delta'$. For the first term of the RHS of the above inequality, we have that
\begin{equation}\label{ineq:diam}
    \begin{split}
        P_{T,x;\mathcal{L}^{}}\left(\eta_k > k^{1+\varepsilon} \right) & \le  P_{T,x;\mathcal{L}^{}}\left(\eta_k > k^{1+\varepsilon}, \mathsf{diam}(T_{\tau_{k-1}})<k^{\varepsilon/2} \right)\\
        & \quad+ \Pd_{T,x;\mathcal{L}^{}}\left(\mathsf{diam}(T_{\tau_{k-1}}) > k^{\varepsilon/2} \right).
    \end{split}
\end{equation}
Observe that, on the event $\{\mathsf{diam}(T_{\tau_{k-1}})<k^{\varepsilon/2}\}$, $\eta_k$ is the optimal strong stationary time of a random walk on a tree with $k-1$ vertices and diameter smaller than $k^{\varepsilon/2}$. By Corollary \ref{cor:mixingbound}, this tree has mixing time at most $k^{1+\varepsilon/2+o(1)}$.
Thus, by the strong Markov property and Proposition~\ref{prop:ssttail}, it follows that there exist two universal constants $C$ and $C'$ such that
\begin{equation}\label{ineq:diam1}
    P_{T,x;\mathcal{L}}\left(\eta_k > k^{1+\varepsilon}, \mathsf{diam}(T_{\tau_{k-1}})<k^{\varepsilon/2} \right) \le Ce^{-C'k^{\varepsilon/3}}.
\end{equation}
Now we are left with the second term of the RHS of \eqref{ineq:diam}. Similarly to~(\ref{eq:taunprob}), but now using that $k^{1+\delta} \gg k$, a Chernoff bound  says that there exists a positive constant $C$ depending on $\gamma$ and $\delta$ such that
\begin{equation}\label{ineq:tautail2}
    \Pd_{T,x;\mathcal{L}}\left(\tau_{k-1} \ge k^{\frac{1+\delta}{1-\gamma}} \right) = \Pd_{T,x;\mathcal{L}}\left(\sum_{i=1}^{k^{\frac{1+\delta}{1-\gamma}}} Z_i \le k-1\right) \le e^{-Ck^{1+\delta}}.
\end{equation}
On the other hand, using that the diameter is nondecreasing, Theorem~\ref{thm:diam} gives us that for large enough $k$ the following holds:
\begin{equation}\label{ineq:diamlarge}
    \begin{split}
         \Pd_{T,x;\mathcal{L}^{}}\left(\mathsf{diam}(T_{\tau_{k-1}}) > k^{\varepsilon/2}, \tau_{k-1} < k^{\frac{1+\delta}{1-\gamma}} \right) & \le \Pd_{T,x;\mathcal{L}^{}}\left(\mathsf{diam}\left(T_{k^{\frac{1+\delta}{1-\gamma}}}\right) > k^{\varepsilon/2}\right)\\
         & \le \frac{1}{k^{2}}.
    \end{split}
\end{equation}
Thus, using \eqref{ineq:tautail2} and \eqref{ineq:diamlarge} we obtain that 
$$
\Pd_{T,x;\mathcal{L}^{}}\left(\mathsf{diam}(T_{\tau_{k-1}}) > k^{\varepsilon/2} \right) \le k^{-2}+e^{-Ck^{1+\delta}},
$$
which combined with \eqref{ineq:first} and \eqref{ineq:diam1} gives us the upper bound
$$
P_{T,x;\mathcal{L}^{}}\left(\eta_k \geq \Delta\tau_k  \right) \le k^{-1-\delta'} + Ce^{-C'k^{\varepsilon/3}}+ k^{-2}+e^{-Ck^{1+\delta}}.
$$
The proof of this case is then finished by the Borel-Cantelli lemma.\\

\noindent {\it \underline{Case $\gamma = 1$.} } The proof of this case follows the same approach as before, but it uses some different bounds due to the fact that now $\sum_{j=1}^k j^{-\gamma}$ is not of polynomial order in $k$. Similarly to the previous case, we start with the following bound
\begin{equation}\label{ineq2:first}
    \begin{split}
        P_{T,x;\mathcal{L}^{}}\left(\eta_k \geq \Delta\tau_k  \right) & \le P_{T,x;\mathcal{L}^{}}\left(\eta_k > k^3 \right) + \Pd_{T,x;\mathcal{L}^{}}\left(\Delta \tau_k \le k^3 \right).
    \end{split}
\end{equation}
Recall that $|T_{\tau_{k-1}}| = k-1$, which implies trivially that $\mathsf{diam}(T_{\tau_{k-1}}) \le k-1$ and that $T_{\tau_{k-1}}$ has a mixing time bounded by $Ck^2$, where $C$ an absolute constant. Then, by Proposition \ref{prop:ssttail} together with Strong Markov property the following upper bound holds 
$$
P_{T,x;\mathcal{L}^{}}\left(\eta_k > k^3 \right) \le e^{-C'k},$$ 
for another absolute constant $C'$. As for the second term of the RHS of \eqref{ineq2:first}, we begin with the following upper bound:
\begin{equation}\label{ineq2:second}
    \Pd_{T,x;\mathcal{L}^{}}\left(\Delta \tau_k \le k^3 \right) \le \Pd_{T,x;\mathcal{L}^{}}\left(\Delta \tau_k \le k^3, \tau_{k-1} \ge k^6 \right) + \Pd_{T,x;\mathcal{L}^{}}\left(\tau_{k-1} \le k^6 \right).
\end{equation}
Now recall that, by the definition of $\tau_k$,
\begin{equation*}
   \Pd_{T,x;\mathcal{L}^{}}\left(\tau_{k-1} \le k^6  \right) = \Pd_{T,x;\mathcal{L}^{}}\left(\sum_{j=1}^{k^6} Z_j \ge k-1 \right),
\end{equation*}
with independent Bernoulli variables $Z_j \sim \mathsf{Ber}(p_j)$.
On the other hand, we have that this sum of Bernoullis has expected value bounded from above by $6 G \log k$. Thus, again by \cite[Theorem 2.7]{mcdiarmid1998concentration}, we obtain that 
\begin{equation}\label{ineq2:sec}
    \Pd_{T,x;\mathcal{L}^{}}\left(\tau_{k-1} \le k^6  \right) \le e^{-Ck},
\end{equation}
for some constant $C$. 

Now, observe that, on the event $\tau_{k-1}\ge k^6$, the increment $\Delta \tau_k$ dominates a geometric distributed random variable of parameter $\G k^{-6}$, because, by the definition of $\cL$, the probability of adding a new leaf is at most $\G k^{-6}$ in each step. This leads us to the following upper bound:
\begin{equation}\label{ineq2:t}
    \Pd_{T,x;\mathcal{L}^{}}\left(\Delta \tau_k \le k^3, \tau_{k-1} \ge k^6 \right) \le \Pd\left(\mathsf{Geom}(\G k^{-6}) \le k^3 \right) \le \frac{C}{k^3},
\end{equation}
for some absolute constant $C$. Finally, putting all the pieces together, one may obtain an upper bound of form
$$
P_{T,x;\mathcal{L}^{}}\left(\eta_k \geq  \Delta\tau_k  \right) \le \frac{C}{k^3},
$$
which, combined with Borel-Cantelli lemma, gives us the result for $\gamma = 1$.
\end{proof}

\subsection{
Application of the Transfer Principle: proof of Corollaries \ref{cor:degdist}-\ref{cor: kth.level}
}
As applications of the Transfer Principle (Theorem \ref{thm:generaltransfer}), we will prove Corollaries \ref{cor:degdist}, \ref{cor:height}, \ref{cor:maxdeg} and \ref{cor: kth.level}.

\begin{remark}[Edge-seed vs.~loop-seed]\label{rem:edge-loop}
The asymptotic properties addressed in the Corollaries are known to hold $P_{G_0}^{\rm edge}$-almost surely for the BA-tree (see reference in each proof).  However, since the hypothesis in Theorem~\ref{thm:generaltransfer} requires  an asymptotic property holding $P_{G_0^{\rm loop}}$-almost surely for the BA-tree, some extra care is needed to assure that, for an asymptotic property $\mathcal{G}$, $P_{G_0^{\rm edge}}\left(\{G_t\}_{t \in \mathbb{N}} \in \mathcal{G}\right)=1 \implies P_{G_0^{\rm loop}}\left(\{G_t\}_{t \in \mathbb{N}} \in \mathcal{G}\right)=1$. 
We will show that this is actually the case for the properties considered in the Corollaries, but each of them requires a specific argument.  However, the common ground to relate   $P_{G_0^{\rm edge}}$-almost surely to $P_{G_0^{\rm loop}}$-almost surely is  the following argument, due to T. F. M\'ori. Intuitively, we ``stretch out'' the vertex with a loop $G_0^{\rm loop}$ into $G_0^{\rm edge}$
and one end of the loop is
colored red, the other is colored blue. 
Then, we color the edges later attached to the root by red or blue according to whether they ``joined the red or the blue end'' of the loop. 
This way a red and a blue tree is created from the root. We do the same with $G_0^{\rm edge}$, but it is even simpler because $G_0^{\rm edge}$ has two vertices and to one of them  ($v_0$) a red tree will be attached, and to the other ($v_1$) a blue tree. 
Hence, given a graph sequence of the BA-tree starting from $G_0^{\rm edge}$, by merging together $v_0$ and $v_1$ (and transforming the initial edge between them in a loop), we assign another graph sequence that will be distributed as the BA-tree starting from $G_0^{\rm loop}$.
\end{remark}

\begin{proof}[Proof of Corollary \ref{cor:degdist} {\rm (Power-law Degree Distribution)}:] For a graph $G$, and $d\geq 1$, let us define 
$$
N_G(d) := \# \text{ of vertices having degree exactly } d \text{ in }G,
$$
and $|G|$ the number of vertices in $G$. Also, write $p_d := \frac{4}{d(d+1)(d+2)}$. Let $\mathcal{G}_{p_d}$ be the following graph property
$$
\mathcal{G}_{p_d} := \left \lbrace \{G_t\}_t \,: \, \lim_{t\to\infty}\frac{N_{G_t}(d)}{|G_t|} = p_d\right \rbrace.
$$
For any fixed $d$,  $\mathcal{G}_{p_d}$ is an asymptotic graph property which
 holds $P_{G_0^{\rm edge}}$-almost surely  for the BA-tree, see, e.g.  \cite{mori2002random}. The proof will follow from Theorem~\ref{thm:generaltransfer} if we show that $\mathcal{G}_{p_d}$ also holds $P_{G_0^ {\rm loop}}$-almost surely. 
 Let $\widehat{N}_{G_t}(d)$ denote the number of vertices having degree exactly $d$ after merging together $v_0$ and $v_1$ (see, Remark~\ref{rem:edge-loop}). Then for any fixed $d$ it holds that 
 \[
N_{G_t}(d) -2 \leq \widehat{N}_{G_t}(d)\leq N_{G_t}(d) +1, 
 \]
 since, after merging $v_0$ and $v_1$,  the number of vertices of degree $d$ may not change, increase at most by one (if $d_{G_t}(v_0)+d_{G_t}(v_1)=d$) or decrease at most by two  (if $d_{G_t}(v_0)=d_{G_t}(v_1)=d$). Dividing then by $|G_t|-1$ (i.e. the size of the graph after contracting the initial edge) and using the coupling discussed in Remark~\ref{rem:edge-loop}, we conclude that $\mathcal{G}_{p_d}$ also holds $P_{G_0^ {\rm loop}}$-almost surely. 
\end{proof}

\begin{proof}[Proof of Corollary \ref{cor:height} {\rm (Height Process)}:] By Theorem~8.1 of \cite{remcoVol2} it follows that 
$$
\mathcal{G} :=\left \lbrace \{G_t\}_t \,: \, \lim_{t\to\infty}\frac{\text{height of } G_t}{\log|G_t|} = \frac{1}{2c} \right \rbrace
$$
is an asymptotic graph property for the BA-tree which holds $P_{G_0^{\rm loop}}$-almost surely, where $c$ is the solution of $ce^{c+1}=1$. Then the proof follows exactly the same lines as  the proof of Corollary \ref{cor:degdist}.
\end{proof}

\begin{proof}[Proof of Corollary \ref{cor:maxdeg} {\rm (Maximum Degree)}:] Let $\mathcal{G}_{\rm max}$ be the following asymptotic graph property
\begin{equation}
    \mathcal{G}_{\rm max} := \left \lbrace \exists \lim_{k\to \infty}\frac{\mathsf{max. deg}(G_k)}{\sqrt{|G_k|}} \text{ a.s. and it is positive.} \right \rbrace
\end{equation}
In \cite{mori2005maximum}, the author shows that $\mathcal{G}_{\rm max}$ holds for the BA-tree $P_{G_0^{\rm edge}}$-almost surely, that is,
\begin{equation}\label{eq:maxdeg}
    \lim_{k\to \infty}\frac{\mathsf{max.deg}(G_k)}{\sqrt{|G_k|}}  = \zeta, \, P_{G_0^{\rm edge}}\text{-a.s.}.
\end{equation}
In order to use Theorem \ref{thm:generaltransfer}, we need to show that $\mathcal{G}_{\rm max}$ holds for the BA-model $P_{G_0^{\rm loop}}$-almost surely.

In order to do so, let $\{G_k^{\rm edge}\}_k$ and $\{G^{\rm loop}_k\}_k$  be the BA-model started from  $G_0^{\rm edge}$ and $G_0^{{\rm loop}}$ respectively. By Remark \ref{rem:edge-loop}, the maximal degree of $G_t^{{\rm loop}}$ can be obtained from $G_t^{\rm edge}$ by contracting the initial edge $G_0^{\rm edge}$. This means that if $v_1$ and $v_2$ are the two initial vertices of $G_0^{\rm edge}$, the following holds
\begin{equation}\label{eq:maxdegloop}
    \mathsf{max.deg}(G_t^{{\rm loop}}) = \max \left \lbrace \mathsf{max.deg}(G_t^{{\rm edge}}), \mathsf{deg}(v_1) + \mathsf{deg}(v_2)\right \rbrace.
\end{equation}
Using \eqref{eq:maxdeg} and Theorem 2.1 in \cite{mori2005maximum}, which implies \eqref{eq:maxdeg} for the degree of a fixed vertex, it follows that there exists a strict positive random variable $\zeta'$ such that 
\begin{equation}
    \lim_{k\to \infty}\frac{\mathsf{max.deg}(G_k)}{\sqrt{|G_k|}}  = \zeta', \, P_{G_0^{\rm loop}}\text{-a.s.}.
\end{equation}
Finally, by Theorem \ref{thm:generaltransfer} the result follows.
\end{proof}
\begin{proof}[Proof of Corollary \ref{cor: kth.level} {\rm (Lower levels of the tree)}:] Let $k\ge 1,\ d\ge 1$, and $q_d=\frac{1}{d(d+1)}$ and consider the following asymptotic property 
\begin{align}\label{eq: Mori.is,smart}
\mathcal{G}_{q_d} :=\left \lbrace \{G_t\}_t \,: \,    \lim_{n\to\infty}\frac{X[n,k,d]}{X[n,k]}=q_d \right \rbrace.
\end{align}
By Theorem 3.1 in \cite{mori2006surprising}, the property $\mathcal{G}_{q_d}$ above 
holds $P_{G_0^ {\rm edge}}$-almost surely for the BA-tree, where the level are computed with respect to a fixed vertex of the initial edge. The proof of this corollary will follow from Theorem~\ref{thm:generaltransfer} if we show that $\mathcal{G}_{q_d}$ also holds $P_{G_0^ {\rm loop}}$-almost surely. 
Starting from $G_0^{\rm edge}$, and having Remark~\ref{rem:edge-loop} in mind, let us assume (without loss of generality)  that the levels are computed with respect to the red vertex of the initial edge, i.e., the red vertex $v_0$ is the unique vertex at level $0$ while the blue one $v_1$ is at level $1$.
For $k\geq 1$ and $d\geq 1$, let  $X[n,k,d]$ be the number of vertices at level $k$ with degree $d$ at time $n$. Also, let  $R[n, k,d]$ ($B[n,k,d]$) denote the number of red (blue) vertices at level $k$ with degree $d$ at time $n$.  
Let $\widehat{X}[n,k,d]$ be the number of vertices at level $k$ with degree $d$ at time $n$ after merging together $v_0$ and $v_1$. Then, for every $k\geq 1$ and $d\geq 1$, it holds that
\[
\widehat{X}[n,k,d] = R[n, k,d] + B[n, k+1,d]. 
\]
What we are after is to show that $\frac{\widehat{X}[n,k,d]}{\sum_{d\geq 1}\widehat{X}[n,k,d]}=q_d$, $P_{G_0^ {\rm loop}}$-almost surely.  
Let  $r_n$ ($b_n$) denote the ``local time of the red (blue) tree'', i.e.,  the number of edges attached to the red (blue) tree by time $n$. An easy application of the Borel-Cantelli lemma shows that, almost surely,
$\lim_{n\to\infty}r_n=\lim_{n\to\infty}b_n=\infty.$
Then, applying \cite[Theorem 3.1]{mori2006surprising} to \eqref{eq: Mori.is,smart} in the red and the blue trees separately, along with the fact that $r_n,b_n\to\infty$ a.s., we obtain that for each fixed $k'\ge 1,d'\ge 1$, 
\begin{align*}
  \lim_{n\to\infty}\frac{R[r_n,k',d]}{\sum_{d\geq 1}R[r_n,k',d]}=q_d, \text{ and }   \lim_{n\to\infty}\frac{B[r_n,k',d]}{\sum_{d\geq 1}B[r_n,k',d]}=q_d, \text{ a.s.}
\end{align*}
Observing that $R[n,k,d]= R[r_n,k,d]$ and $B[n,k+1,d]= B[r_n,k,d]$, for each fixed $k\ge 1,d\ge 1$, we conclude that 
\begin{align*}
    &\frac{\widehat{X}[n,k,d]}{\sum_{d\ge 1}\widehat{X}[n,k,d]}=\frac{R[r_n,k,d]+B[b_n,k,d]}{\sum_{d\ge 1}\widehat{X}[n,k,d]}\\
   & =\frac{R[r_n,k,d]}{\sum_{d\ge 1}R[r_n,k,d]}\cdot\frac{\sum_{d\ge 1} R[r_n,k,d]}{\sum_{d\ge 1}\widehat{X}[n,k,d]}+\frac{B[b_n,k,d]}{\sum_{d\ge 1}B[b_n,k,d]}\cdot\frac{\sum_{d\ge 1} B[b_n,k,d]}{\sum_{d\ge 1}\widehat{X}[n,k,d]}\\
   & =:\frac{R[r_n,k,d]}{\sum_{d\ge 1}R[r_n,k,d]}\gamma_{n,k,d}+\frac{B[b_n,k,d]}{\sum_{d\ge 0}B[b_n,k,d]}(1-\gamma_{n,k,d})\to q_d,
\end{align*}
$P_{G_0^ {\rm loop}}$-almost surely. 
\end{proof}

\section{Proof of Theorem \ref{thm:diam}: Upper bound on the diameter}\label{sec:diameter}

We are left with the proof of Theorem \ref{thm:diam} in order to finally prove Theorem \ref{thm:transfer}. This will be done in this section, after recalling some basic notions about graphs.

\begin{definition}[Order, size, diameter]
We recall that the \emph{order}  of a graph  is the cardinality of its vertex set, while the \emph{size} of a graph is the cardinality of its edge set. The \emph{diameter} of the graph is the maximum of all graph distances between pairs of vertices in the graph.
\end{definition}

Before stating the main result of the section, we warn the reader that constants denoted by capital letters may change from lemma (proposition, etc.) to lemma.

The proof of Theorem~\ref{thm:diam} is quite involved. It follows from the lemma below, whose proof we postpone to Subsections~\ref{subsec: proof.doublediam} and~\ref{sec:step_ii}.

\begin{lemma}[Small time diameter increment]
\label{lem:doublediam} 
Let $\gamma \in (1/2,1]$ and $\delta<2\gamma -1$. 
Fix $T$ (finite initial tree),  $d\ge 3$ (factor), $j$ (time factor) and $n$ sufficiently large.  Then, there exists a $C_1>0$, depending on $\gamma$, $\g$, $\G$, $d$, and $|T|$, and a $\kappa>0$ depending only on $\gamma$ and $\delta$, not on $d$, such that 
\begin{align*}
    \mathbb{P}_{T,x;\mathcal{L}^{}} \left(\mathsf{diam}\left(T_{n/2^{j-1}}\right) - \mathsf{diam}(T_{n/2^j})\ge d\left(\frac{n}{2^j}\right)^{1-2(1-\gamma)-\delta}\right)\le 
    C_1\left(\frac{2^{j}}{n}\right)^{\delta\gamma+\kappa(d/2-1)}\!\!,
\end{align*}
holds for all $n\ge 1$ and $j=1,2,...,\lfloor\log_2 n\rfloor$.
    
Note that if $\kappa(d/2-1)>\gamma-\delta$  (which can be achieved by taking $d$ large enough), then, for any fixed $j\ge 1$, the right-hand side tends to zero as $n\to\infty$.
\end{lemma}
For now, we will focus on showing how Theorem \ref{thm:diam} follows from the above lemma.

\begin{proof}[Proof of Theorem \ref{thm:diam}] Given $\varepsilon>0$ and a large enough $n$, we have that
\begin{equation}\label{eq:diamTn}
    \mathsf{diam}(T_n) = \sum_{j=1}^{(1-\varepsilon)\log_2(n)}\left[\mathsf{diam}(T_{n/2^{j-1}}) - \mathsf{diam}(T_{n/2^{j}})\right] + \mathsf{diam}(T_{n^{\varepsilon}}).
\end{equation}
Now, for a fixed $d\in \mathbb{N}$ and $\delta < 2\gamma -1$, let $A^{(d)}_j$, for $j \in \{1,2,\dots, (1-\varepsilon)\log_2(n)\}$, be the following event:
\begin{equation}
    A^{(d)}_j := \left\lbrace \mathsf{diam}(T_{n/2^{j-1}}) - \mathsf{diam}(T_{n/2^{j}}) \ge d\left(\frac{n}{2^j}\right)^{1-2(1-\gamma)-\delta}\right \rbrace.
\end{equation}
By Lemma \ref{lem:doublediam} and the union bound, there exist a constant $\kappa>0$ depending on $\gamma,\g,\G, \delta$, and a constant $C_1$ depending on $d, \gamma, \g,\G, |T|$, such that, if $d$ is chosen large enough so that $\kappa (d/2-1)>\gamma-\delta$, then 
\begin{equation}\label{eq:unionAj}
    \mathbb{P}_{T,x;\mathcal{L}^{}}\left( \bigcup_{j=1}^{(1-\varepsilon)\log_2n}A_j^{(d)}\right) \le \frac{C_1n^{\gamma-\delta}}{n^{\kappa(d/2-1)}}\sum_{j=1}^{(1-\varepsilon)\log_2n}\frac{2^{j \kappa(d/2-1)}}{2^{j(\gamma-\delta)}} \le \frac{C_2}{n^{\varepsilon[\kappa(d/2-1)-\gamma+\delta]}},
\end{equation}
where $C_2$ depends on  $d, \gamma, \g, \G, |T|$ and $\delta$.  

On the other hand,  notice that the number of vertices in $T_{n^{\varepsilon}}$, hence also its diameter, are at most $n^{\varepsilon}+|T|$.
So, using~\eqref{eq:diamTn} and~\eqref{eq:unionAj},
it follows that, with probability at least~$1-C_2{n^{-\varepsilon[\kappa(d/2-1)-\gamma+\delta]}}$ , 
\begin{equation}
    \mathsf{diam}(T_n) \le  n^{\varepsilon} + |T| + d\sum_{j=1}^{(1-\varepsilon)\log_2n}\left(\frac{n}{2^j}\right)^{1-2(1-\gamma)-\delta} \le  n^{\varepsilon} + |T|+ C_3 dn^{2\gamma-1-\delta},
\end{equation}
for a positive constant $C_3$ depending on $\gamma, \g, \G$ and $\delta$ only, since $\sum_j2^{-j(2\gamma-1-\delta)} < \infty$, using that $\delta<2\gamma-1$.

Now, given $\gamma$, and the $\varepsilon,M>0$ of the statement of Theorem~\ref{thm:diam}, pick a $0<\delta=\delta(\gamma,\veps)< 2\gamma-1$ such that $2\gamma-1-\delta < \veps$. 
Then choose $d=d(\gamma, \delta, M, \varepsilon)$ large enough so that
$\varepsilon[\kappa(\gamma, \delta)(d/2-1)-\gamma+\delta]>M$. With these choices, we define 
$$C:= 1+dC_3(\gamma,\g,\G,\delta).$$ 
This $C$ works for large enough $n$'s and can obviously be modified to work for all $n\ge 1$.
\end{proof}

\subsection{Proof of Lemma \ref{lem:doublediam}}\label{subsec: proof.doublediam}

The rationale behind the proof of Lemma \ref{lem:doublediam} is as follows. Since adding new vertices to the tree becomes more expensive as time goes to infinity ($p_n$ is decreasing), one may expect that in a small time window, the walker cannot increase the diameter too much as it would imply adding many new vertices in a short period of time. 

Formally, what we will show is that, using the proper scale for the time window, the walker is not capable of growing new sub-trees having a large amount of vertices. Consequently, it is unable to increase the diameter by large amounts. This idea is summarized in the next lemma. Again, we will postpone its proof to the next subsection, showing here only how Lemma \ref{lem:doublediam} follows from it.
\begin{lemma}
\label{lem:constantincrdiam}
Let $\gamma \in (1/2,1]$ and $\delta<2\gamma -1$. 
Given $T$ (finite initial tree) and   $d\ge 3$ (diameter increment), there exists a  constant~$C_2>0$ depending on $\gamma, \g, \G$, $d$ and $|T|$,  and there exists a $\kappa>0$  (depending on $\gamma$, $\delta$, but not on $d$), such that for all large $n$'s,
    \[
    \mathbb{P}_{T,x;\mathcal{L}^{}} \left(\mathsf{diam}\left(T_{n+n^{2(1-\gamma)+\delta}}\right) - \mathsf{diam}(T_n)\ge d\right)
    \le \frac{C_2 n^{1-\gamma}}{n^{\kappa (d/2-1)}}.
    \]
\end{lemma}

\begin{proof}[Proof of Lemma \ref{lem:doublediam}] For fixed $n,j \in \mathbb{N}$ and $k \in \{1,2,\dots, (n/2^{j})^{2\gamma-1-\delta}\}$, define
$$\hat n(j,k):={n/2^{j}+(k-1)(n/2^j)^{2(1-\gamma)+\delta}},$$
then let
\begin{equation}\label{eq: long}
    \Delta \mathsf{diam}_{n,j,k} := \mathsf{diam}(T_{\hat n(j,k+1)}) - \mathsf{diam}(T_{\hat n(j,k)})\,.
\end{equation}
Thus,
\begin{equation}
    \mathsf{diam}(T_{n/2^{j-1}}) - \mathsf{diam}(T_{n/2^{j}}) = \sum_{k=1}^{(n/2^{j})^{2\gamma-1-\delta}}\Delta \mathsf{diam}_{n,j,k}. 
\end{equation}
From this it follows that if 
$$
    \mathsf{diam}(T_{n/2^{j-1}}) - \mathsf{diam}(T_{n/2^{j}}) \ge d\left(\frac{n}{2^j}\right)^{1-2(1-\gamma)-\delta},
$$
then there exists $k \in \{1,2, \dots, (n/2^{j})^{2\gamma-1-\delta}\}$ such that
$$
\Delta \mathsf{diam}_{n,j,k} \ge d.
$$
To finish the proof, we use  Lemma \ref{lem:constantincrdiam} and the union bound as follows.
Note that, for $\hat n=\hat n(j,k)$, since the diameter is nondecreasing in time,
$$ \Delta\mathsf{diam}_{n,j,k} \le 
\mathsf{diam}(T_{\hat n+\hat n^{2(1-\gamma)+\delta}}) - \mathsf{diam}(T_{\hat n}).
$$
Hence, using that $\hat{n} < n/2^{j-1}$,
\begin{align*}
&\mathbb{P}_{T,x;\mathcal{L}^{}} \left( \Delta \mathsf{diam}_{n,j,k} \ge d\right)\\
&\ \ \ \le \mathbb{P}_{T,x;\mathcal{L}^{}} \left(\mathsf{diam}(T_{\hat n+\hat n^{2(1-\gamma)+\delta}}) - \mathsf{diam}(T_{\hat n})\ge d\right)\stackrel{\text{Lemma}\ \ref{lem:constantincrdiam}}{\le}C_2 \hat n^{1-\gamma-\kappa (d/2-1)}\\
&\ \ \ < C_2  (2^{j-1}/n )^{\gamma-1+\kappa (d/2-1)}.
\end{align*}
Finally, the union bound leads to the bound
    \begin{equation}
    \begin{split}
    \label{leads.to}
        (n/2^{j})^{2\gamma-1-\delta}&\cdot C_2 \, (2^{j-1}/n )^{\gamma-1+\kappa (d/2-1)}= \frac{C_2}{2^{\gamma - 1 + \kappa(d/2-1)}}  (2^{j}/n )^{\delta-\gamma+\kappa (d/2-1)}
        \\
        &\leq C_1 \, (2^{j}/n )^{\delta-\gamma+\kappa (d/2-1)},
        \end{split}
    \end{equation}
with $C_1=C_2 2^{1-\gamma}$ that depends only on $\gamma, \g, \G$, $d$ and $|T|$, as required.
%
\end{proof}

\subsection{Upper bound on the order of sub-trees: Proof of Lemma \ref{lem:constantincrdiam}} \label{sec:step_ii}

To complete the proof of Theorem \ref{thm:diam} we are left with proving Lemma \ref{lem:constantincrdiam}. The proof is quite involved and depends on controlling the order of sub-trees created by vertices added in a certain time window, which is done in Proposition \ref{lem:step_ii}.

Recall that all trees generated by the \name\ process have a common root. Thus, given a vertex $v \in T_t$ for some time $t\ge 0$, we have a natural notion of the descendants of $v$.

\begin{definition}[Descendants]
The descendants of the vertex $v$ are formed by those vertices that are further away from the root than $v$, and such that the only path connecting them to the root contains $v$. 
\end{definition}

With this idea in mind, one can  define {\it the tree of the descendants of $v$ at time $t$}. 

\begin{definition}[The tree of the descendants of $v$ at a given time]
Let $t\ge 0$. The tree $T_t(v)$ is the sub-tree of $T_t$ formed by the vertex $v$ and all its descendants. 
\end{definition}

We will use the notation $T(v)$ when the time at which we refer to the sub-tree of descendants of $v$ is irrelevant.

For $k \in \mathbb{N}$ we let $v_k$ be the vertex possibly added at time $k$, depending on the value of $Z_k$ (if $Z_k=0$, the vertex $v_k$ is not added to the tree). 
Let  $T_t(v_k)$ be  the sub-tree of descendants of $v_k$ at time $t$ and  let $D_{k,t}$ denote the order of $T_t(v_k)$; that is, 
\begin{equation}\label{eq:deg}
    D_{k,t} := \begin{cases}
        |T_t(v_k)|, & \text{ if }Z_k = 1; \\
        0, & \text{otherwise}.
    \end{cases}
\end{equation}
If $t<k$  we set $T_t(v_k)=\emptyset$ and if  vertex $v_k$ is not added (i.e., $Z_k=0$),  we set $T_t(v_k)=\emptyset$  for all $t$. 

For a fixed $n$, $T_n$ will be called the ``blue tree'' and every sub-tree created after $n$ will be called a ``red sub-tree.''
The next result states that in a time window of length $n^{2(1-\gamma) +\delta}$, the order  of any particular red sub-tree is small.

\begin{proposition}
\label{lem:step_ii}
Let $\gamma \in (1/2,1]$ and $\delta<2\gamma -1$. 
Fix  $T$ (finite initial tree) and  $d\in \mathbb N$ (order).
Then, there exist $C_1=C_1(\gamma,\g,\G,d,|T|)$ and $\kappa=\kappa(\gamma,\delta)>0$  ($\kappa$ does not depend on $d$, and $C_1$ does not depend on $\delta$),   such that for all large $n$'s:
\[
\mathbb{P}_{T,x;\mathcal{L}^{}} (D_{k,n+ n^{\myexp}} \ge d)\le \frac{C_1}{k^\gamma n^{\kappa (d-1)}}, \  k= n,n+1,\dots\;.
\]
\end{proposition}

Before we prove Proposition \ref{lem:step_ii}, we will show how it implies Lemma \ref{lem:constantincrdiam}. 

\begin{proof}[Proof of Lemma \ref{lem:constantincrdiam}] 
Introduce the shorthand $N(n,\gamma,\delta):=n+n^{\myexp}$.
Notice that if, for a given $\delta < 2\gamma-1$ and $d \in \mathbb{N}$, we have that
$
\mathsf{diam}\left(T_{N(n,\gamma,\delta)}\right) - \mathsf{diam}(T_n)\ge d,
$
then there must be a red sub-tree with diameter at least $d/2$. To see this, consider the shortest path between two vertices of $T_{N(n,\gamma,\delta)}$ that realizes $\mathsf{diam}\left(T_{N(n,\gamma,\delta)}\right)$. It must have at least $d$ red edges. Since the red edges come from at most two red sub-trees, one of these sub-trees must have diameter at least $d/2$. Thus, by Proposition \ref{lem:step_ii} and the union bound, it follows that, for some positive constant $C_2$ depending on $\gamma, \g, \G$, $|T|$  and $d$ only,
\begin{align*}
    \mathbb{P}_{T,x;\mathcal{L}^{}} \left(\mathsf{diam}\left(T_{N(n,\gamma,\delta)}\right) - \mathsf{diam}(T_n)\ge d\right) 
& \leq
\mathbb{P}_{T,x;\mathcal{L}^{}} \left(\bigcup_{k=n}^{N(n,\gamma,\delta)}\left \lbrace D_{k,M(n,\gamma,\delta)} \ge d/2 \right \rbrace \right)\\
& \le \sum_{k=n}^{N(n,\gamma,\delta)}\frac{C_1}{k^\gamma n^{\kappa (d/2-1)}} \le \frac{C_2 n^{1-\gamma}}{n^{\kappa (d/2-1)}},
\end{align*}
 using that $N(n,\gamma,\delta) < 2n$. This proves Lemma~\ref{lem:constantincrdiam}.
\end{proof}

\subsubsection{Proof of Proposition~\ref{lem:step_ii}}

The proof of  Proposition~\ref{lem:step_ii} relies on three  auxiliary results, namely Corollary~\ref{lem:fixedtree} (a consequence of \cite[Lemma 3.3]{englander2021tree}) and  Lemmas ~\ref{lemma:expecnst} and ~\ref{lemma:rec},  which will be presented below. 
 
The first statement is  a general result for simple random walk on {\em fixed} trees which will be useful in the sequel.

\begin{corollary}[Bound for fixed trees]\label{lem:fixedtree} Let $T$ be a finite rooted tree, $v$ a vertex of $T$ (it may also be the root) and $T(v)$ the sub-tree of descendants of $v$. For $t\in\mathbb N$, let  $N_t$ be the number of visits to $T(v)$ in $t$ steps of a simple random walk on~$T$. Then, 
\[
E_v\left[N_t\right] \le |T(v)|\left[\frac{t+3}{|T|-1}+48|T|\right]\;
\]
\end{corollary}

\begin{proof} Let $N_t(u)$ be the number of visits to a given vertex $u$ of $T(v)$ in $t$ steps, so that
$$N_t = \sum_{u \in T(v)}N_t(u).$$
By Lemma 3.3 of \cite{englander2021tree},
$$
E_v\left[N_t(v)\right] \le \mathsf{deg}_T(v)\left[\frac{t+3}{2(|T|-1)}+24|T|\right].
$$ 
Furthermore, for any vertex $u$ and any time $t$, we have $E_v\left[N_t(u)\right] \leq 
E_u\left[N_t(u)\right]+1$, since, after the first visit to $u$, we can apply the strong Markov property. Thus,
$$
E_v\left[N_t(u)\right] \le \mathsf{deg}_T(u)\left[\frac{t+3}{2(|T|-1)}+24|T|\right]+1.
$$
Summing this up over $u\in T(v)$, and noticing that $T(v)$ being  a subtree of descendants implies that 
$$
\sum_{u \in T(v)} \mathsf{deg}_T(u) \le 2(|T(v)|-1),
$$
give us 
$$
E_v\left[N_t\right] \le \big(|T(v)|-1\big)\left[\frac{t+3}{|T|-1}+48|T|\right]+|T(v)|.
$$
This implies the claimed bound using $|T(v)|\le |T|$.
\end{proof}

 Since in the \name, the domain of the walk changes during the walk, we will need to keep track of the random times when new vertices are added and the times when  their sub-trees increased. For this reason we will need some  further  definitions.
 \begin{definition}[First time reaching order $d$]Let $k\ge 1$ be given.
For $d \in \mathbb{N}$, define inductively the following sequence of stopping times:
\begin{equation*}
    \eta_{k,1} := \begin{cases}
        k, & \text{ if }Z_k = 1; \\
        +\infty, & \text{ otherwise}\;,
    \end{cases}
\end{equation*}
and for $d\ge 2$,
\[
\eta_{k,d} :=  \begin{cases}
\inf \{ t > \eta_{k,d-1} \; : \; D_{k,t} = d\}\;, & \text{ if $\eta_{k,d-1} <\infty$;
}
\\
+\infty\;, & \text{ otherwise}\;.
\end{cases}
\]
 So, $\eta_{k,d}$ is the first time $T(v_k)$ reaches order  $d$ (if $v_k$ is not added then $\eta_{k,d}=+\infty$, for all $d\geq 1$). 
 \end{definition}
 
Given $k\in \mathbb{N}$ (vertex index), $t,s \in \mathbb{N}$ (times) with $k\geq t$   and $d\geq 1$ (order),   we denote by  $N_{k,t,t+s}^{(d)}$  the number of visits to $T(v_k)$ in the time interval $[t,t+s]$ at times when $T(v_k)$ has order exactly $d$. Formally, we make the following definition.
\begin{definition}[Number of visits to subtree of $v_k$]Given $k\in \mathbb{N}$, $t,s \in \mathbb{N}$ with $k\geq t$   and $d\geq 1$, let   

\begin{equation}\label{eq:N}
    N_{k,t,t+s}^{(d)} := \begin{cases}
        \sum_{j=t}^{t+s}\mathbb{1}\{X_j \in T_j(v_k), D_{k,j}=d\}\;, & \text{ if }Z_k = 1;\\
        0\;, & \text{ otherwise}\;.
    \end{cases}
\end{equation}
\end{definition}
In the definition above we require $k\geq t$ since we are interested in the number of visits to sub-trees  of vertices created after time $t$.  
For convenience, when $d=1$, we drop the superscript $(d)$ in the above definition.
%

With respect to $N_{k,t,t+s}^{(d)}$, we have the following result.  
\begin{lemma}[Bound on mean number of visits to subtree]\label{lemma:expecnst} Fix $T$ (finite initial tree), $k$ (vertex index), $d$ (order), and $t$ (time) with $ k\ge t$. Then there exist positive constants $C_1,C_2,C_3$ and $C_4$, depending on $\gamma,\g,\G$ only, such that, for any $s\ge 1$,
 \begin{equation*}
        \begin{split}
            \mathbb{E}_{T,x;\mathcal{L}^{}} \left[N_{k,t,t+s}^{(d)}\right] &\le C_1d\left[\frac{s}{|T|+C_2t^{1-\gamma} -1}+d+|T| + (t+s)^{1-\gamma}\right]\mathbb{P}_{T,x;\mathcal{L}^{}}\left( \eta_{k,d} < t+s\right) \\
            & \quad + C_4d(|T| + t+s) e^{-C_3t^{1-\gamma}}\,.
        \end{split}
    \end{equation*}
\end{lemma}

\begin{proof} First, by  \eqref{eq:N}, we have the identity $N_{k,t,t+s}^{(d)} =\mathbb{1}\{\eta_{k,d} < t+s\} N_{k,t,t+s}^{(d)}$. Moreover,
    \begin{equation*}
      \mathbb{1}\{\eta_{k,d} < t+s\} N_{k,t,t+s}^{(d)} \le \mathbb{1}\{\eta_{k,d} < t+s\} N_{k,\eta_{k,d},\eta_{k,d}+s}^{(d)}\;, \quad \mathbb{P}_{T,x;\mathcal{L}^{}}\text{-a.s.}
    \end{equation*}
    (This is because on the right-hand side we count the visits in a time window of size exactly $s$, whereas on the left-hand side the time window has at most size $s$.) Thus, by the strong Markov property it follows that 
    \begin{equation}\label{eq:smarkovnst}
      \mathbb{E}_{T,x;\mathcal{L}^{}} \left[N_{k,t,t+s}^{(d)}\right] \le \mathbb{E}_{T,x;\mathcal{L}^{}}\left[ \mathbb{E}_{T_{\eta_{k,d}}, X_{\eta_{k,d}}; \mathcal{L}^{(\eta_{k,d})}} \left[N_{k,0,s}^{(d)}\right]\mathbb{1}\{\eta_{k,d} < t+s\} \right]\,,
    \end{equation}
where it is important to point out that $v_k$ in $N_{k,0,s}^{(d)}$ is not the $k$-th vertex added by the shifted process, but the original vertex $v_k$ that belongs to $T_{\eta_{k,d}}$ on the event $\{\eta_{k,d} < t+s\}$.
   
To bound the expression
    \begin{equation}\label{handle}
    (*):=\mathbb{E}_{T_{\eta_{k,d}}, X_{\eta_{k,d}}; \mathcal{L}^{(\eta_{k,d})}} \left[N_{k,0,s}^{(d)}\right] \mathbb{1}\{\eta_{k,d} < t+s\},
    \end{equation} 
    we use coupling.
    To this end, 
    given $T_{\eta_{k,d}}$ and $X_{\eta_{k,d}}$, let $(W,P_{X_{\eta_{k,d}}})$ denote a simple random walk on $T_{\eta_{k,d}}$, starting from $X_{\eta_{k,d}}$, and  let $N^{W}_{s}$ be the number of visits to the sub-tree $T_{\eta_{k,d}}(v_k)$ by the walker $W$ in $s$ steps.
     By ignoring some parts of $X$, we are going to couple  
    $\{X_i\}_{0\le i\le s}$  with $\{W_i\}_{0\le i\le s}$    (setting $X_0=W_0=X_{\eta_{k,d}}$) in such a way that  all possible subsequent visits of $X$ to $T_{\eta_{k,d}}(v_k)$ counted by $N_{k,0,s}^{(d)}$ are counted by $N^{W}_{s}$ as well.
   To do this, notice that: 
   \\
   (i) $T_{\eta_{k,d}}(v_k)$ may reach order $d+1$ before completing $s$ steps, and hence the count in $N_{k,0,s}^{(d)}$ stops,
   \\
   (ii) $X$ may leave $T_{\eta_{k,d}}$ at other parts of the tree for ``excursions'' on the larger tree. 
   \\
   In  case (ii) we just delete those excursions outside of $T_{\eta_{k,d}}$, while in  case (ii) we must stop the walk $W$. Consequently, the walk $(W,P_{X_{\eta_{k,d}}})$  on $T_{\eta_{k,d}}$ obtained this way,   may take {\it less} than $s$ steps. To make up for that, and have exactly $s$ steps for $W$,  simply ``complete'' the remaining steps of  $W$ by letting it walk on $T_{\eta_{k,d}}$, independently of $X$. 
    It then follows that, given $T_{\eta_{k,d}}$ and $X_{\eta_{k,d}}$,
   one has
    \begin{equation}\label{itfollowsthat}
      (*)\le E_{X_{\eta_{k,d}}} \left[ N^{W}_{s} \right] \, \mathbb{1}\{\eta_{k,d} < t+s\}.
    \end{equation}
    This is of course also true without the indicator, but it is weaker without that, and we will later use the stronger version.
    Now, in order to bound the expected value in \eqref{itfollowsthat}, we distinguish between two cases, according to  whether $X_{\eta_{k,d}} \in T_{\eta_{k,d}}(v_k)$ or not. (Note that $X_{\eta_{k,d}} \notin T_{\eta_{k,d}}(v_k)$ is indeed possible. This is because the process might add a vertex and {\it then} jump to the father of $v_k$. Recall that first we change the tree at time $n$, and only then we move the walker to get $X_n$). In the first case, using the more general notation $(W,P_{v}), v\in T_{\eta_{k,d}}$, and letting $H^W_{v_k}$ be the hitting time of $W$ to $v_k$, by the strong Markov property of  $W$, using the commute time identity \cite[Proposition 10.7]{levin2017markov} for the first term on hitting times on trees, and using Corollary \ref{lem:fixedtree} for the second term, it follows that, given $T_{\eta_{k,d}}$ and $X_{\eta_{k,d}}$, on the event $\{\eta_{k,d} \le t+s\}$, we have
    \begin{equation}
        \begin{split}
            E_{X_{\eta_{k,d}}} \left[ N^{W}_{s} \right] \le E_{X_{\eta_{k,d}}} \left[ H^W_{v_k} \right] + E_{v_k} \left[ N^{W}_{s} \right] 
            & \le 2d^2 + d\left(\frac{s+3}{|T_{\eta_{k,d}}|-1}+ 48|T_{\eta_{k,d}}|\right)\\
            & \le C_1\left[d^2 + d\left(\frac{s}{|T_{\eta_{k,d}}|-1}+ |T_{\eta_{k,d}}|\right)\right],
        \end{split}
    \end{equation}
    with an absolute constant $C_1<\infty$. If $X_{\eta_{k,d}} \notin T_{\eta_{k,d}}(v_k)$, then the above bound follows without the term $d^2$.
    Thus, on the event $\{\eta_{k,d} \le t+s\}$, we can continue the inequality \eqref{itfollowsthat} as
\begin{equation*}
    \begin{split}
        (*) \le E_{X_{\eta_{k,d}}} \left[ N^{W}_{s} \right] 
        &\le C_1\left[d^2 + d\left(\frac{s}{|T_{\eta_{k,d}}|-1}+ |T_{\eta_{k,d}}|\right)\right] \\
        &  \stackrel{k\geq t}{\le} C_1\left[d^2 + d\left(\frac{s}{|T_{t}|-1}+ |T_{t+s}|\right)\right] \,,
    \end{split}
    \end{equation*}
since the assumption $k\ge t$ implies that $T_t\subset T_{\eta_{k,d}}$, and on the event $\{\eta_{k,d} \le t+s\}$ we have that $T_{\eta_{k,d}} \subset T_{t+s}$. Plugging this into \eqref{eq:smarkovnst} leads to 
    \begin{equation}\label{eq:expts}
      \mathbb{E}_{T,x;\mathcal{L}^{}} \left[N_{k,t,t+s}^{(d)}\right] \le \mathbb{E}_{T,x;\mathcal{L}^{}}\left[ C_1\left[d^2 + d\left(\frac{s}{|T_{t}|-1}+ |T_{t+s}|\right)\right]\mathbb{1}\{\eta_{k,d} \le t+s\}\right].
    \end{equation}
     Given the initial condition $(T,x)$, since the tree growth is governed by the  sequence of laws $\mathcal{L}^{}=\{\mathsf{Ber}{p_n}\}_{n\ge 1}$ satisfying Definition~\ref{def:G}, it follows that 
    \[
    |T_t| =  |T| + \sum_{r=1}^tZ_{r} \implies  \mathbb{E}_{T,x;\mathcal{L}^{}} \left[ |T_t| \right] = 
    |T| + \Theta \left(t^{1-\gamma}\right),
    \]
    where the constants involved in the $\Theta$ notation depend only on $\gamma,\g,\G$. Next, using Chernoff bounds for the independent indicators $Z_i$, one sees that there exist some positive constants $C_2, C'_2, C_3$ and $C'_3$ depending on $\gamma,\g,\G$ only such that 
    \begin{align*}
    \mathbb{P}_{T,x;\mathcal{L}^{}} \left( |T_t| \le |T| + C_2 t^{1-\gamma}\right) &\le e^{-C_3t^{1-\gamma}};\\ \mathbb{P}_{T,x;\mathcal{L}^{}} \left( |T_{t+s}| \ge |T| + C'_2(t+s)^{1-\gamma}\right) &\le e^{-C'_3t^{1-\gamma}}.
    \end{align*}
    Combining these inequalities with \eqref{eq:expts}, one obtains
    \begin{equation*}
        \begin{split}
            \mathbb{E}_{T,x;\mathcal{L}^{}} \left[N_{k,t,t+s}^{(d)}\right] &\le C_1d\left[\frac{s}{|T|+C_2t^{1-\gamma} -1}+d+|T| + C'_2 (t+s)^{1-\gamma}\right]\mathbb{P}_{T,x;\mathcal{L}^{}}\left( \eta_{k,d} < t+s\right) \\
            & \quad + C_4d(|T| + t+s) e^{-C''_3t^{1-\gamma}},
        \end{split}
    \end{equation*}
where all constants depend on $\gamma,\g,\G$ only.
  After arranging the constants if necessary, the proof of the lemma is complete.
\end{proof}

The next result is a recursion for the tail probability of  $D_{k,t+s}$.

\begin{lemma}[Recursive tail bound for  $D_{k,t+s}$]\label{lemma:rec} 
 Fix $T$ (finite initial tree),  $k$ (vertex index), $d$ (order), $t,s \ge 1$ (times) with $ k\ge t$ and fix  $\kappa < \gamma$. Then, there exist some positive constants $C_1$, $C_2$, $C_3$ and $C_4$ depending on $\gamma,\g,\G$ only  (same as in Lemma~\ref{lemma:expecnst}), such that 
$$
\mathbb{P}_{T,x;\mathcal{L}^{}} (D_{k,t+s} \ge d+1)\le I+II\,,
$$ 
where
\begin{align*}
 & I:= \left \lbrace \frac{\G\, (t+s)^{\gamma}}{\g\, t^{\gamma}}\left[1 - \left( 1- \frac{\g}{(t+s)^\gamma}\right)^{t^{\gamma-\kappa}}\right]
       +\Delta\right \rbrace \mathbb{P}_{T,x;\mathcal{L}^{}}\left( D_{k,t+s} \ge d\right)
\end{align*}
with
\begin{equation*}
    \Delta := C_1d\left[\frac{s}{|T|+C_2 t^{1-\gamma} -1}+d+|T| + (t+s)^{1-\gamma}\right]t^{-\gamma+\kappa}\,,
\end{equation*}
           and
           $$ II:= C_4d(|T| + t+s) e^{-C_3t^{1-\gamma}}t^{-\gamma + \kappa}\,.
           $$
\end{lemma}

\begin{proof} Let $F^{k,d}_{t,s}:=\{N_{k,t,t+s}^{(d)} \le t^{\gamma - \kappa}\}.$ Clearly, 
$
  E^{k,d}_{t+s}:=  \left \lbrace \eta_{k,d} \le t + s \right \rbrace = \left \lbrace D_{k,t+s} \ge d \right \rbrace, 
$
and 
    \begin{equation}\label{eq:dktged}
        \mathbb{P}_{T,x;\mathcal{L}^{}} (E^{k,d+1}_{t+s}) 
         \le \mathbb{P}_{T,x;\mathcal{L}^{}} (E^{k,d+1}_{t+s}\cap  F^{k,d}_{t,s})+ \mathbb{P}_{T,x;\mathcal{L}^{}} ( [F^{k,d}_{t,s}]^c)\;.
    \end{equation}
   Markov's inequality along with Lemma~\ref{lemma:expecnst} yield 
    \begin{align}
    \label{eq:b1}
     \mathbb{P}_{T,x;\mathcal{L}^{}} ( [F^{k,d}_{t,s}]^c)
  \le \Delta\mathbb{P}_{T,x;\mathcal{L}^{}}\left( \eta_{k,d} < t+s\right) + II.
    \end{align} 
    For the first term of the RHS of \eqref{eq:dktged}, fix $j \le t^{\gamma-\kappa}$. 
    Let $H'_i$ denote the time of the $i$-th visit to $T(v_k)$ after time $\eta_{k,d}$, i.e., after reaching order $d$. We then have the  identity
    \begin{align*}
    \mathsf{Event}_j:&=E^{k,d+1}_{t+s}\cap  \lbrace N_{k,t,t+s}^{(d)} = j\rbrace \\  
    &= \left \lbrace  Z_{H'_j + 1} = 1, Z_{H'_{j-1}+1} = 0, \dots, Z_{H'_1+1} = 0, \eta_{k,d} < t+s \right \rbrace.
    \end{align*}
    In words, $\mathsf{Event}_j$ occurs if $T(v_k)$ has order at least $d+1$ before time $t+s$ and has been visited $j$ times while having order $d$. Hence, at each of these $j$ steps on $T(v_k)$, in the next step the order has failed  to increase  to $d+1$, and only succeeded after the $j$-th step, that is, at time $H'_j +1$. The failures and the success are described formally by the $Z_{H'_i+1}$'s. Since the $Z$'s are independent and $p_n$ is decreasing in $n$, 
    \begin{equation*}
      \mathbb{P}_{T,x;\mathcal{L}^{}} (\mathsf{Event}_j )
     \le \frac{\G}{t^\gamma}\left( 1 - \frac{\g}{(t+s)^{\gamma}}\right)^{j-1}\mathbb{P}_{T,x;\mathcal{L}^{}}\left( E^{k,d}_{t+s}\right).
    \end{equation*}
    Summing over $j$ from $1$ to $t^{\gamma - \kappa}$ leads to 
    \begin{align*}
       & \mathbb{P}_{T,x;\mathcal{L}^{}} (E^{k,d+1}_{t+s}\cap \{N_{k,t,t+s}^{(d)} \le t^{\gamma - \kappa}\}) \\
        &\qquad \qquad\le \frac{\G\, (t+s)^{\gamma}}{\g\, t^{\gamma}}\left[1 - \left( 1- \frac{\g}{(t+s)^\gamma}\right)^{t^{\gamma-\kappa}}\right] \mathbb{P}_{T,x;\mathcal{L}^{}}\left( E^{k,d}_{t+s}\right),
    \end{align*}
which, combined with \eqref{eq:b1}, proves the result.
\end{proof}

We are finally  ready to prove Proposition~\ref{lem:step_ii}.

\begin{proof}[Proof of Proposition~\ref{lem:step_ii}]
Setting $t:= n$ and $s:=n^{\myexp}$, with $\delta<2\gamma -1$  and $\kappa<\gamma$ in Lemma~\ref{lemma:rec}, and  using  the shorthand  $D^{(\delta)}_{n,k}:=D_{k,n+n^{2(1-\gamma)+\delta}}$ and $u(d):=\mathbb{P}_{T,x;\mathcal{L}^{}} (D^{(\delta)}_{n,k} \ge d)$, we obtain from Lemma~\ref{lemma:rec} that  
    \begin{equation}
       \begin{split}
       \label{eq:recursive_U}
       u(d+1)\leq&
      \frac{\G}{\g}\left(\frac{ n+ n^{\myexp} }{n}\right)^\gamma\left\{1 - \left(1-\frac{\g}{(n+n^{\myexp})^\gamma}\right)^{n^{\gamma -\kappa}} \right\} u(d) 
      \\
      & + \frac{C_1d}{n^{\gamma-\kappa}}\left[\frac{n^{\myexp}}{C_2 n^{1-\gamma}}+d+|T| + (n+ n^{\myexp})^{1-\gamma}\right] u(d) 
      \\
      &+ C_4d(|T| + n+ n^{\myexp}) e^{-C_3n^{1-\gamma}}n^{-\gamma + \kappa},
       \end{split}
    \end{equation}
    where the positive constants $C_1,C_2,C_3$ and $C_4$, depend on $\gamma,\g,\G$ only.
Since $2(1-\gamma) +\delta < 1$, we have that $n+ n^{\myexp}\leq 2n$. 
Furthermore, 
\begin{equation*}
    \left(1-\frac{\g}{(n+n^{\myexp})^\gamma}\right)^{n^{\gamma -\kappa}} \geq    \left(1-\frac{\g}{n^\gamma}\right)^{n^{\gamma -\kappa}}  \ge e^{-\frac{3\g}{2} n^{-\kappa }} \geq 1-\frac{3\g}{2} n^{-\kappa }, 
\end{equation*}
for all $n\ge 1$ if $\g$ is small enough, which of course we may assume. Consequently,
\begin{equation*}
    \frac{\G}{\g} \left(\frac{ n+ n^{\myexp}}{n}\right)^\gamma\left\{1 - \left(1-\frac{\g}{(n+n^{\myexp})^\gamma}\right)^{n^{\gamma -\kappa }} \right\} \le \frac{3\G}{2^{1-\gamma}n^\kappa} =: \frac{C'_1}{n^\kappa},
\end{equation*}
where $C_1'$ depends on $\gamma,\G$ only. In light of the above, from~\eqref{eq:recursive_U} we obtain that  
\begin{align*}
u(d+1)\leq  \left(\frac{ C'_1}{n^{\kappa}}+C_2'\frac{dn^{\myexp}}{n^{1-\kappa}} 
+ \frac{C''_1 d^2}{n^{\gamma-\kappa }} + \frac{C_1d (2n)^{1-\gamma}}{n^{\gamma-\kappa}} \right) u(d)
       + \frac{ 3 C_4 d n}{n^{\gamma - \kappa}}e^{-C_3n^{1-\gamma}},
\end{align*}
where $C''_1= C_1 (|T|+1)$  is a positive constant depending on $\gamma,\g,\G$ and $|T|$ and $C_2'=C_1/C_2$ depends on $\gamma,\g,\G$ only. In the last term, by assuming  $n$ sufficiently large, we used that $|T|\leq n$. 
Now, since $\gamma \in (1/2,1]$ and $\delta < 2\gamma-1$, we can choose $\kappa=\kappa(\gamma, \delta)$ sufficiently small such that 
$1-\kappa - [2(1-\gamma)+\delta]>\kappa$, and hence $\gamma - \kappa -(1-\gamma)>\kappa$ and $\gamma - \kappa >\kappa$ all happen, and thus
we can further simplify~\eqref{eq:recursive_U} to obtain that, for all sufficiently large $n$,
\begin{align*}
u(d+1)\leq  \frac{ \widehat{C}_1 d^2}{n^{\kappa}} u(d)
       + 3C_4 d n^{1-\gamma + \kappa}e^{-C_3n^{1-\gamma}}\leq\frac{ \widehat{C}_1 d^2}{n^{\kappa}} u(d)
       + d e^{-\widehat{C}_2 n^{1-\gamma}},
\end{align*}
   where $\widehat{C}_1$ depends on $\gamma,\g,\G$ and $|T|$ and $\widehat{C}_2$ is a positive constant depending on $\gamma,\g,\G$ only.  

Note that, from the definition of $u(d)$, along with the  assumption that  $2(1-\gamma)+\delta<1$,  we have that $u(d)=0$ unless $d\le 2n$  and $k\le n+n^{\myexp}$. Hence, we  will only consider $d$'s (for a given $n$) that do not exceed $2n$. 
It follows that, for all $n$ sufficiently large,
\begin{align}\label{eq:rec}
  u(d+1)\le  \widehat{C}_1 d^2n^{-\kappa} u(d)
       + e^{-Ln^{1-\gamma}},
\end{align}
where $L:=\widehat{C}_2/2$. 
Continuing the recursion \eqref{eq:rec} in $d$ backwards, all the way back to $u(1) =k^{-\gamma}$, we obtain that, for all $n$ sufficiently large, 
\begin{align*}
     u(d+1)&\le   \frac{1}{k^\gamma}\prod_{j=1}^{d} 
     \left(\widehat{C}_1 j^2n^{-\kappa}\right) + e^{-Ln^{1-\gamma}}\sum_{j=0}^{d-1}
     \left(\widehat{C}_1 d^2\right)^j n^{-\kappa j}
     \\
      &\leq \frac{1}{k^\gamma} 
     \left(\widehat{C}_1 d^2\right)^d n^{-d\kappa} + \widehat{C}_3e^{-Ln^{1-\gamma}},
\end{align*}
where $\widehat{C}_3:= \sum_{j=0}^{d-1}\left(\widehat{C}_1 d²\right)^j$ is a constant depending on $\gamma, \g,\G, |T|$ and $d$. Moreover, since $k\le n+n^{\myexp}$, by replacing $\left(\widehat{C}_1d^2\right)^{d}$ by a larger constant if necessary, there exists a positive constant $C'$ depending on $\gamma, \g,\G, |T|$ and $d$ only, such that for all $k \in [n, n+n^{\myexp}]$ and $n\ge 1$,
$$
\widehat{C}_3\exp\left\{-Ln^{1-\gamma}\right\} \le \frac{1}{k^\gamma} 
     C' n^{-d\kappa},
$$
which implies that 
\begin{equation*}
    u(d+1) \le \frac{1}{k^\gamma} 
     C'n^{-d\varepsilon} + \widehat{C}_3\exp\left\{-Ln^{1-\gamma}\right\} \le  \frac{2C'}{k^\gamma} n^{-d\kappa},
\end{equation*}
completing the proof.
\end{proof}


%

\section*{Acknowledgments }
J. E. is indebted to Y. Peres for first bringing the model to his attention and for further helpful conversations, and to T. F. M\'ori for answering  some questions concerning his work on the BA model.

J. E.'s research was partially supported by Simons Foundation Grant 579110. G.P.~was supported by the ERC Consolidator Grant 772466 ``NOISE'' until January 2024 to find the main ideas, by the ERC Synergy Grant 810115 ``DYNASNET'' from February 2024 to find applications, and also by the Hungarian National Research, Development and Innovation Office, OTKA grant K143468, throughout. G.I.~was supported by FAPERJ (grant E-26/210.516/2024).

\bibliographystyle{amsplain}
\bibliography{ref.bib}

 


\end{document}